\documentclass[final,5p,times]{elsarticle}

\usepackage{amssymb}
\usepackage[T1]{fontenc}
\usepackage[utf8]{inputenc}
\usepackage{microtype}
\usepackage{graphicx}
\graphicspath{{ElsArticleTemplate/}}
\usepackage[outdir=./]{epstopdf}
\usepackage{amsmath}
\usepackage{mathtools}
\usepackage{siunitx}
\usepackage{longtable,tabularx}
\usepackage{multicol}
\usepackage[export]{adjustbox}
\usepackage{float}
\usepackage{multirow}
\usepackage[table]{xcolor}

\usepackage{lineno}
\modulolinenumbers[5]

\usepackage{caption}
\usepackage{subcaption}
\usepackage{url}
\usepackage{float}
\usepackage{svg}

\usepackage{algorithm}
\usepackage{algpseudocode}
\usepackage{tikz}
\usetikzlibrary{shapes,arrows,positioning,fit,backgrounds,calc}
\usepackage{pgfplots}
\usepackage{import}
\pgfplotsset{compat=1.16}
\usepackage{hyperref}
\hypersetup{pdfauthor={Name}}
\usepackage{tcolorbox}
\usepackage{xr}
\usepackage{booktabs} 
\usepackage{stfloats}
\usepackage{makecell}

\usepackage{xcolor}  
\definecolor{myblue}{RGB}{18, 159, 255}
\definecolor{mypurple}{RGB}{125, 46, 143}
\definecolor{mygreen}{RGB}{120, 171, 49}
\usepackage[normalem]{ulem}  
\colorlet{add}{blue}

\graphicspath{{img/}}
\DeclareSymbolFontAlphabet{\mathbb}{AMSb}

\journal{Acta Astronautica}

\begin{document}

\begin{frontmatter}



\title{A Methodology for Integrating Life Cycle Assessment into a Multidisciplinary Design Analysis and Optimization Framework for Sustainable Launcher Development}


\author[inst1,inst3]{Alice De Oliveira\corref{cor1}}
\ead{alice.de-oliveira@isae-supaero.fr}
\cortext[cor1]{Corresponding Author}
\author[inst2]{Mathieu Balesdent}
\ead{mathieu.balesdent@onera.fr}
\author[inst2]{Lo\"{\i}c Brevault}
\ead{loic.brevault@onera.fr}
\author[inst1]{Annafederica Urbano}
\ead{annafederica.urbano@isae-supaero.fr}

\affiliation[inst1]{organization={F\'{e}d\'{e}ration ENAC ISAE-SUPAERO ONERA, Universit\'{e} de Toulouse},
            city={Toulouse},
            country={France}}

\affiliation[inst2]{organization={ONERA/DTIS, Universit\'{e} Paris-Saclay},
            city={Palaiseau},
            country={France}}

\affiliation[inst3]{organization={Centre national d’études spatiales (CNES)},
            city={Paris},
            country={France}}

\begin{abstract}
The increasing number of orbital and sub-orbital launches makes it necessary to investigate the environmental impacts of launch vehicles and to incorporate eco-design considerations into their development. In response to these emerging concerns, the European Space Agency has promoted the use of  the Life Cycle Assessment (LCA) methodology as a standardization to support the mitigation of environmental impacts associated with present and future space missions. 
This need is further amplified in the NewSpace, where start-ups and established industrial actors explore numerous configurations and innovative technologies, reinforcing the importance of integrating environmental considerations.
At early design stages, defining the launch vehicle architecture may be formalized through the solving of a multi-physics optimization problem based on Multidisciplinary Design Analysis and Optimization (MDAO) methods. 
The different physics involved are numerically modeled through disciplines (e.g., propulsion, aerodynamics, structure, and trajectory) that are coupled together inside an MDAO process to obtain trade-offs of candidate configurations according to different performance criteria.
In this paper, a methodology is proposed to integrate an LCA discipline within an MDAO framework dedicated to launch vehicle design. The approach relies on the definition of parametric life-cycle inventories that depend on design and coupling variables modeling the interactions between the disciplines. It includes the production of the stages components and propellants, as well as the associated transport to the launch site. Furthermore, launch emissions are evaluated from the optimized trajectory profiles and characterized in terms of climate change impact.
The approach is illustrated on the design of a representative expendable launch vehicle showing the capabilities of the proposed MDAO framework enhanced with predictive LCA. Specifically, multi-objective optimizations are performed to assess the tradeoffs between a launch vehicle performance criterion and environmental impact indicators. The results highlight the antagonistic behaviors among environmental impact categories, emphasizing the importance of carefully defining environmental objectives when conducting eco-design studies for launch vehicles. Overall, the generic nature of the methodology lays the first foundation for integrating LCA into launch vehicle early-stage design processes, thereby enabling exploration of trade-offs between performance, cost, and environmental considerations.
\end{abstract}




\end{frontmatter}


\section{Introduction}\label{sec:intro}
The space sector is undergoing profound industrial and technological changes. Reusable launch systems have fundamentally altered the economics of space access, reducing launch costs and reshaping mission design and planning. This shift has contributed to an emerging of a more competitive and diversified space economy, involving both governmental agencies and a growing number of influential private actors pursuing more autonomous and cost-efficient access to orbit. In parallel, commercial space applications have expanded rapidly, most notably through deployment of mega-constellations providing global broadband internet coverage. While these developments are often associated with improved efficiency, they also contribute to a rebound effect, where reduced launch costs lead to increased launch frequency and, consequently, a growing overall environmental burden. In fact, recent analyses suggest that the climate impact of space activities could increase by a factor of nine by 2050 \cite{Miraux2022}. 

Despite increasing scientific attention to space debris, the broader environmental footprint of launch vehicles, whether reusable or expendable, remains insufficiently quantified and requires further investigation \cite{DominguezCalabuig2024}. Beyond climate change, relevant impacts include ozone depletion, resource scarcity (metals, minerals, and fossil resources), and land and water use. While adjacent sectors such as commercial aviation are actively integrating sustainability considerations into system design \cite{keiser2023lca}, the space industry still lacks comprehensive regulatory frameworks addressing environmental performance. 

In response to these emerging concerns, the European Space Agency (ESA) has promoted the use of the Life Cycle Assessment (LCA) methodology as a standardization to support the mitigation of environmental impacts associated with present and future space missions \cite{esa2025lca}. The LCA framework has become an increasingly relevant tool in the aerospace sector, enabling the quantification of environmental burdens across all stages of a product’s life cycle, from raw material extraction to end-of-life management \cite{maury2020application}. Among existing approaches, the Environmental Footprint (EF) method developed by the European Commission \cite{Damiani2022} offers a harmonized framework within the European context. 
Nevertheless, most available LCA studies are carried out after the completion of vehicle design \cite{gallice2018environmental, Harris2019, Pellegrino2025}, which prevents them from guiding environmentally informed decisions during the conceptual and preliminary design phases.
In contrast, NewSpace private companies are beginning to explore sustainable launch vehicle configurations. For example, the European company MaiaSpace has indicated that it is assessing the environmental impact of its future mini launcher using LCA methodology from an early design stage \cite{Vila2025}.

In fact, the early design of a launch vehicle is inherently a multidisciplinary process involving coupled disciplines such as aerodynamics, structural design, propulsion, and trajectory optimization \cite{Castellini2014, Balesdent2022}. The resulting Multidisciplinary Design Analysis and Optimization (MDAO) framework typically aims to minimize an objective function, such as a mass or cost criterion (e.g., Gross Lift-Off Weight, GLOW), while satisfying multiple constraints (e.g., maximum axial load, maximum dynamic pressure). 
Recent studies have begun investigating how to integrate environmental assessments within MDAO frameworks to identify designs that balance performance and sustainability \cite{bellier:hal-04188708, TORMENA2022, Musso2024, Gregorio2026}. 
The main limitation of these approaches lies in their only partial integration of the LCA methodology, both in terms of the life-cycle phases considered and the range of impact categories analyzed. Indeed, these studies primarily focus on climate change impacts, assessed through radiative forcing and the associated Global Warming Potential (GWP) indicator, while being restricted to the launch flight phase. As a result, they overlook the environmental impacts associated with other stages, notably the production and assembly of the launch vehicle. 
This highlights the need for a comprehensive MDAO framework for launch vehicle design that fully integrates LCA-based environmental considerations.
However, embedding an LCA discipline within such an optimization process poses challenges related to model fidelity (i.e., the level of detail and accuracy of the LCA models), data availability, and compatibility with other physical and engineering models. Traditional LCA requires detailed information on the definition of the launch vehicle and therefore on all its life-cycle phases (e.g., R\&D, manufacturing, transport, etc.), data that are typically unavailable at early design stages. 
Consequently, as for other traditional disciplines involved in early-stage design, the level of modeling and the hypotheses adopted for LCA must be adapted to remain meaningful in this context.
Further limitations arise from the limited availability of launch-vehicle-specific datasets (e.g., materials, propellants, manufacturing processes, end-of-life strategy). Moreover, modeling emissions across multiple atmospheric layers remains complex, and the environmental impacts of high-altitude emissions and space operations are largely unquantified \cite{Ross2014}. 

These challenges underscore the need for a structured approach to incorporate LCA at early design stages.
This paper builds upon the methodology introduced by the same authors in De Oliveira et al.~\cite{DeOliveira2025} to integrate LCA as a discipline within an MDAO framework and therefore enable informed trade-offs between vehicle performance, cost, and environmental impacts. It relies on the definition of parametric life-cycle inventories, which depend on the design and coupling variables of the MDAO framework and enable fast computation of the LCA impacts. More particularly, it includes the production of the stages' components and propellants from raw materials with the transport phase associated (e.g., production of an engine at a manufacturing site then transported to the launch site in a different location). The result is therefore a predictive LCA discipline suitable for MDAO integration. 
In this paper, the methodology is enhanced to account for corresponding atmospheric impacts from launch emissions thanks to existing models in the literature that enable their computation \cite{FischerLEATegusphere-2026-3050}, and the available metrics to characterize their impacts on the climate change \cite{IPCC2014, Lee2021, Sherwood2018, Hauglustaine2022, Lammel1995}. The capabilities of this integrated predictive LCA discipline into a launch vehicle design proces are illustrated for a medium-lift Two-Stage-To-Orbit (TSTO) LOx/LCH\textsubscript{4} (Liquid Oxygen / Liquid Methane) vehicle, extending its application from the sensitivity and parametric analyses of \cite{DeOliveira2025} to a multi-objective optimization setting.
Therefore, the main contribution of this paper is the proposal of a mature methodology for integrating LCA as a discipline within an MDAO framework, enabling eco-design from a holistic, multidisciplinary perspective rather than on a discipline-by-discipline basis. This integration also enables the obtained framework to support not only optimization \cite{Balesdent2025}, but also other early-stage studies, including sensitivity and parametric analyses, as well as uncertainty propagation \cite{Brevault2025}.

The paper is organized as follows. Sec.~\ref{sec:2} introduces the LCA methodology and specifies the actual state-of-the-art of the environmental considerations for space activities. More particularly, the main challenges and limitations are highlighted in the context of early design stages. Then, Sec.~\ref{sec:3} describes the methodology developed in this paper to integrate a predictive LCA discipline within an early-stage MDAO framework for launch vehicles. While the MDAO context is described, a focus is made on the development of the LCA discipline with the development of the parametric life-cycle inventories, the estimation of the launch emissions, and the subsequent evaluation of the corresponding environmental impacts. In order to illustrate the capabilities of the methodology developed in this paper, the LCA discipline is embedded into a generic MDAO framework in Sec.~\ref{sec:4}, and the resulting optimization problem is formulated as a multi-objective problem minimizing the GLOW and an environmental impact indicator.

\section{LCA methodology and current practices for launch vehicles}
\label{sec:2}
This section introduces the LCA methodology and examines how it is currently applied to space activities. It provides an overview of the standard practices, identifies the main challenges encountered in modeling space systems at early design stages, and highlights the knowledge gaps that remain in quantifying environmental impacts, particularly for high-altitude emissions.

\subsection{LCA methodology}
\label{sec:2-lca-methodology}
\begin{figure}[!h]
    \centering
    \includegraphics[width=\linewidth]{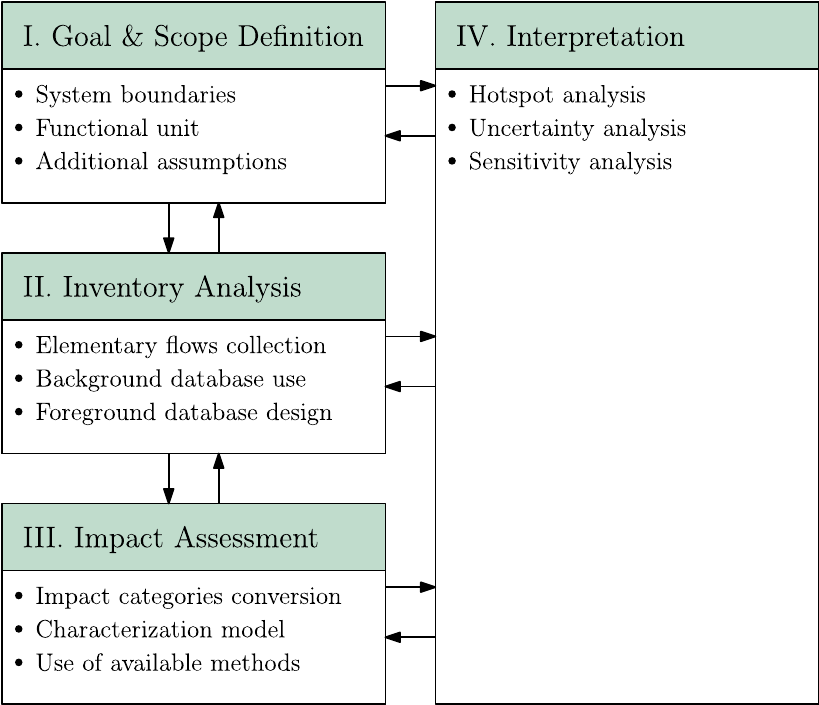}
    \caption{Description of the phases of the LCA methodology, adapted from ISO 14040:2006 \cite{iso14040}. The double arrows suggest an iterative process.}
    \label{fig:lca-metodo}
\end{figure}

The LCA methodology is grounded in two international standards: ISO 14040:2006, which provides the principles and framework \cite{iso14040}, and ISO 14044:2006, which specifies requirements and guidelines \cite{iso14044}. It is based on four phases described in Fig.~\ref{fig:lca-metodo}: (i) goal and scope definition, (ii) inventory analysis, (iii) impact assessment, and (iv) interpretation. 

The first phase defines the boundaries of the system to be assessed. It should therefore mention whether an activity is included in the assessment of the launch vehicle impacts, such as for example the ground segment infrastructure or the R\&D phase. Additionally, a so-called functional unit must be defined. The latter corresponds to a quantification of the system's function as well as the definition of associated assumptions with respect to the scenario studied. It must serve as a reference unit for comparison of several systems, here launch vehicles. 
Maury et al.~\cite{maury2020application} provide a comprehensive overview of the studies and developments related to the application of the LCA method in the space sector. Much of the work in this field originates from Europe, and particularly from the ESA through its Clean Space initiative, launched in 2012 \cite{esa_cleanspace_blog}. As part of this effort, ESA developed the ESA LCA Handbook \cite{esa2025lca}, which establishes methodological guidelines for conducting LCA specifically for space systems. Within this document, the system boundaries of space-related activities are carefully defined and examples of functional units are given. The one of a launch vehicle is given by \cite{esa2025lca}:
\begin{description}
\centering
    \item[] \it To transport a [payload of V kg] from destination A to destination B in T time.
\end{description}
where A represents the launch site coordinates, B the nature of the orbit into which the payload is transported to, V is the maximum payload mass that the vehicle can place into orbit, and T is the time during which or at which the service is provided.
Note that the choice of the functional unit and the assumptions associated is a critical stage that must be carefully investigated. 

Then, the second phase consists of the inventory of all the processes that occur throughout the life cycle of the system under study and the collection of the corresponding data. More specifically, the various physical flows including resources, emissions, energy, waste, and materials are quantified for each process. The output is the Life-Cycle Inventory (LCI), which lists all the elementary flows (e.g., CO$_2$, crude material) that are associated with the functional unit. 
This step is particularly resource- and time-consuming, which is why LCI databases exist (e.g., ecoinvent \cite{frischknecht2005ecoinvent}), with aggregated information for unit processes that have their entire life cycle already modeled (e.g., the production of 1 kg of aluminum AL2219). One can note that most existing databases are established for non-space activities. More particularly, it is common practice to combine \textit{foreground} data, i.e., processes based on in-house primary information about the system under study (such as specific materials used or distances traveled), with \textit{background} data which consist of generic unit processes provided by available LCI databases (such as energy mixes, transportation modes, or external product manufacturing).
ESA has regularly updated its space LCI database \cite{ESA_LCA}, first released in 2017. This database compiles space-specific data on materials, advanced manufacturing processes, and propellants. It was initially derived from the ecoinvent database \cite{frischknecht2005ecoinvent}, developed by the Swiss Centre for Life-Cycle Inventories, but has since been extended to include processes unique to space applications. Building on these initiatives to develop dedicated space-related datasets, Wilson \cite{wilson2019thesis} also created the Strathclyde Space Systems Database (SSSD). This resource integrates data from several background databases, such as ecoinvent and the ESA LCA database, to provide ready-to-use LCI for complete space systems, including launch vehicles and satellites. 

The third phase defined by the LCA methodology is the Life-Cycle Impact Assessment (LCIA). It involves converting all generated elementary flows (i.e., the emissions released to and the resources extracted from the environment) into environmental impact categories. Each impact category (e.g., "Climate change") uses a characterization model with a defined indicator (e.g., the GWP over 100 years, GWP100). However, to date, this step is not standardized and several LCIA methods are available in the literature, using their own impact categories and corresponding characterization models, such as ReCiPe \cite{ReCiPe-Huijbregts2016} or IMPACT World+ \cite{IMPACT-Bulle2019}. For space-related activities, the European Commission has engaged the Product Environmental Footprint (PEF) initiative associated with the EF method \cite{Damiani2022} as the standard method to be used in Europe. This method is also recommended by the ESA LCA Handbook \cite{esa2025lca}. In fact, the European Commission is currently working on a specific "Product Environmental Footprint Category Rules for the space sector (PEFCR4space)" with a consortium of private space companies, international and national organizations, non-governmental organizations, and academic partners\footnote{Product Environmental Footprint Category Rules (PEFCR) for the space sector. \href{https://defence-industry-space.ec.europa.eu/eu-space/product-environmental-footprint-category-rules-pefcr-space-sector_en}{\texttt{https://defence-industry-space.ec.europa.eu/eu-space/}}}. 

Finally the fourth phase is the interpretation of the results, the assessment of their robustness to uncertainties along all the previous phases, and the identification of the potential areas for improvement of the system. 

\subsection{Coupling LCA with MDAO for launch vehicles}
\label{sec:2-lca-mdao}
A critical challenge in space-related LCA is the integration of environmental considerations at the earliest stages of system design. In most existing work, LCA is applied once the vehicle architecture has already been established, thereby limiting its impact on early-stage design choices \cite{gallice2018environmental, Harris2019, Pellegrino2025}. By contrast, Miraux et al.~\cite{Miraux2022} underlined the importance of performing environmental assessments in parallel with technical and economic evaluations during the conceptual design of space systems. This integration can be achieved within a MDAO framework for conceptual launcher design studies. While numerous examples of MDAO for launch vehicle design exist \cite{Balesdent2011, Wilken2024}, they typically do not incorporate LCA. In contrast, integrated frameworks including environmental modules are more established in aeronautics. For instance, Pollet et al.~\cite{Pollet2023, Pollet2025} embedded an LCA-based environmental module into a preliminary design framework for unmanned aerial vehicles and aircraft, enabling environmentally informed design decisions. For launch systems, this approach has been initiated by Dominguez Calabuig et al.~\cite{DominguezCalabuig2024, DominguezCalabuig2022}, who assessed the environmental impacts of fleets of launch vehicles representing potential future designs. These vehicles were defined for identical missions and determined using an optimal staging methodology \cite{DominguezCalabuig2022b}. This comprehensive analysis provides a promising foundation for incorporating environmental considerations into early-stage launcher design through the inclusion of an LCA discipline and a suitably formulated optimization problem.

The integration of environmental considerations within MDAO frameworks for space applications remains largely underexplored, with only a limited number of studies addressing this aspect to date. Tormena et al.~\cite{TORMENA2022} investigated the inclusion of radiative forcing \cite{Ross2014} as an LCA-related discipline within an MDAO framework and compared launcher design options using different propellant combinations. Design decisions were then guided by treating the resulting radiative forcing as a constraint in the launch vehicle optimization. More recently, Musso et al.~\cite{Musso2024} developed a complete multidisciplinary framework for sounding rockets. In this approach, the LCA discipline captures the environmental impacts of both fuel production and combustion during launch, while a multi-objective optimization simultaneously maximizes apogee altitude and minimizes climate change impacts, quantified via the GWP100.
Despite the only partial integration of the LCA methodology since the authors do not consider the production and assembly phases of the launch vehicle components, this framework provides a meaningful methodological contribution. Similarly, Gregorio et al.~\cite{Gregorio2026} integrated launch emissions into an early-stage conceptual design framework to assess their impact and guide sustainability-driven design choices; however, they have not explicitly included them in the optimization process. Overall, all these contributions only analyzed the climate change among all the impact categories \cite{Damiani2022} and focused on the propellant use, disregarding the impacts of the production and assembly phases.

Additionally, embedding an LCA discipline within an MDAO framework for early-stage launch vehicle design and defining adequately the optimization problem to be solved is not straightforward \cite{bellier:hal-04188708}. Incorporating environmental metrics requires careful consideration of how the MDAO problem is structured, including the selection of design and coupling variables and the fidelity of the LCA model relative to the other typically lower-fidelity models used at early design stages. Depending on the design objectives, the problem may then be posed as a multi-objective optimization (e.g., performance or cost versus environmental impact) \cite{Musso2024}, a constrained optimization (e.g., limiting environmental impacts below specified thresholds) \cite{TORMENA2022}, or a quality-diversity approach aimed at exploring a diverse set of sustainable design alternatives \cite{Baraton2025}.

As mentioned above, the fidelity of the LCA model must be aligned with the level of available data and knowledge in early design phases. Although LCI databases exist \cite{ESA_LCA, frischknecht2005ecoinvent, wilson2019thesis}, accessing and combining relevant processes with the information obtained from the other disciplines remains resource- and time-intensive, limiting their integration into MDAO frameworks.
Furthermore, for space-related activities, several knowledge gaps still exist, such as the lack of inventory data and the modeling of launch flight emissions as a function of altitude, and therefore the associated environmental impacts in terms of climate change or ozone depletion. 

Another important consideration is the classification of environmental impacts that can be assessed using LCA. In the LCA studies for launch vehicles mentioned above \cite{TORMENA2022,Musso2024,Gregorio2026}, only the climate change impact via the GWP indicator is assessed but many other impact categories exist. In fact, making design decisions based on multiple impact categories can be challenging, as some impacts may be antagonistic. To support decision-making and promote standardization, the European Commission has introduced a normalization and weighting procedure, yielding a single environmental score known as the PEF \cite{sala2017weighting}. According to this methodology, midpoint impact categories such as climate change, water depletion, resources scarcity and land use typically receive the highest weights. Nevertheless, this normalization and weighting approach may not be fully applicable to space-related activities, which are unique in terms of system design, production, operational use, and end-of-life management. For example, Verkammen et al. \cite{Verkammen2023eucass} argued that ozone depletion and resource scarcity should receive greater weighting in space-based LCA based on a European-focused survey of launch activity experts. These findings highlight the need for a more nuanced understanding of single-score methodologies when applied to the environmental assessment of space systems. In fact, the MDAO framework can enable performing optimizations considering multiple environmental indicators (instead of a single score), allowing a more comprehensive analysis.

\subsection{Current projects and research activities}
\label{sec:2-lca-sota}
Several initiatives and working groups have recently been established to advance the application of LCA to space systems, with a particular emphasis on the atmospheric impacts of space activities, which remain insufficiently characterized. Notably, the REACT project funded by the ESA led to the development of the Assessment and Comparison Tool (ACT) \cite{Udriot2023react, udriot2024react} which is a simplified, space-specific, prospective LCA calculation tool designed for early-stage use by engineers, enabling rapid screening of potential environmental impacts for planned missions or space transportation vehicles. 
A key novelty of the tool is its ability to provide preliminary high-altitude atmospheric emissions based on trajectory information. These estimations are adapted from the Launch Emission Assessment Tool (LEAT) and Re-entry Emission Assessment Tool (REAT) developed by the University of Stuttgart \cite{Fischer2025}.
ACT is now maintained and made available to industry by the recently founded spinoff company EcoDeltaV\footnote{EcoDeltaV website. \url{https://ecodeltav.com/}}, which focuses on providing the tool as a service, engineering support, and education in space sustainability.
Additionally, a white paper published by the “LCA of Space Transportation Systems” working group at the University of Stuttgart brings together experts from academia, industry, government, and national agencies to recommend measures for improving understanding of the environmental impacts of launchers \cite{Fischer2024a}.

In fact, to date, launch emissions resulting from propellant combustion are often neglected in the LCA for launch vehicles \cite{gallice2018environmental}, as their contribution to climate change remains poorly quantified.
Nevertheless, a growing number of recent studies \cite{Fischer2025, desain2014potential, Pradon2023, schabedoth2020lca, Ryan2022, Brown2024, Barker2024, Herberhold2026, james2021commercial} have attempted to quantify these impacts more rigorously. Several of these works have compiled emission inventories over various time periods: Desain and Brady \cite{desain2014potential} for 1985-2013, Pradon et al.~\cite{Pradon2023} for 2009-2018, Schabedoth \cite{schabedoth2020lca}, Ryan et al.~\cite{Ryan2022} and Brown et al.~\cite{Brown2024} for the year 2019, Barker et al.~\cite{Barker2024} for 2020-2022, Fischer et al.~\cite{Fischer2025} for 2019-2024, and more recently Herberhold et al.~\cite{Herberhold2026} for the year 2024. These inventories provide first-order estimates of launch vehicle atmospheric emissions, in a manner analogous to inventories developed for the aviation sector \cite{Quadros2022, Teoh2024}.
In addition to these global inventories, some studies \cite{Fischer2025, desain2014potential, Pradon2023, Brown2024, Barker2024, Herberhold2026} have attempted to assign emissions to specific atmospheric layers, producing altitude-dependent emission profiles. Such data are necessary inputs for sophisticated global atmospheric and climate modeling efforts aimed at evaluating impacts on climate change and stratospheric ozone depletion. Running these models generally involves computationally intensive numerical simulations, demanding substantial computing resources.

For practical software implementation, a common methodological framework has emerged to estimate the launch emissions \cite{FischerLEATegusphere-2026-3050, james2021commercial} for early design consideration, based on three interconnected components: (i) estimating primary emission indices, corresponding to the exhaust products released directly from the engine nozzle exit, (ii) evaluating secondary emissions resulting from afterburning chemical reactions with ambient atmospheric constituents, and (iii) determining the vehicle’s trajectory profile (altitude, velocity, and geodetic coordinates) using flight dynamics models. 
Once these emission profiles are estimated, they must be converted into environmental impacts such as climate change or ozone depletion. This step generally requires coupling with the computationally intensive global climate models mentioned above. However, due to their high cost, these models are often replaced by simplified analytical approaches based on characterization factors, which convert emissions into impact indicators such as the GWP100 metric \cite{IPCC2014, Lee2021, Sherwood2018, Hauglustaine2022, Lammel1995}.
Note that this approach assumes that climatic impacts can be estimated from secondary emission inventories alone, implicitly neglecting the intermediate-scale transport and chemical transformation processes that occur between the near-field post-combustion regime and the final atmospheric concentrations relevant at the global climate scale \cite{Fuglestvedt2003}.
Although this general framework is widely adopted in the literature, these estimates remain preliminary and have not been systematically validated, even against higher-fidelity modeling approaches. As a result, the resulting estimates can vary significantly, mainly due to differences in modeling strategies and underlying expert assumptions. Furthermore, most existing studies rely on deterministic approaches, which do not systematically account for the substantial uncertainties inherent to these complex processes.

These advances further underline the importance of reliable characterization approaches capable of converting emissions into quantified environmental impacts, as mentioned above.
Existing LCA studies have typically relied on climate metrics derived from ground- or aviation-based applications \cite{Miraux2022,DominguezCalabuig2022}, with limited consideration of altitude-dependent effects. In particular, the role of emitted species such as aluminum oxide (\(\mathrm{Al_2O_3}\)), black carbon (BC), and gaseous chlorine (\(\mathrm{HCl}\), \(\mathrm{Cl}\), \(\mathrm{Cl_2}\)) remains insufficiently characterized, especially regarding their contributions to climate change and ozone depletion, despite evidence suggesting potentially significant impacts \cite{Ross2014, Ryan2022, Maloney2022}. 
The climatic effects of water emissions at high altitude are also highly uncertain and strongly dependent on altitude \cite{Pletzer:watervapour}, further highlighting the challenges and limitations associated with applying conventional climate metrics to rocket emissions.
In this context, Miraux et al.~\cite{Miraux2022} highlighted key knowledge gaps and explored possible trajectories of space-related environmental impacts up to 2050 under different scenarios. 
More recently, Dominguez Calabuig et al.~\cite{DominguezCalabuig2024, DominguezCalabuig2022} extended this perspective by incorporating updated emission characterization models into a fleet-level LCA of launch vehicles, covering multiple impact categories such as climate change, ozone depletion, water use, and land use. 
Their results show that rocket exhaust emissions can dominate the overall environmental footprint by several orders of magnitude when altitude-dependent characterization factors and improved modeling of BC effects are included.
However, the magnitude of these impacts should be interpreted with caution, as the method used to estimate the updated characterization factors has not been validated in the literature and therefore remains subject to considerable uncertainty.
Overall, the studies discussed in this section highlight both the potential significance of altitude-dependent effects and the considerable uncertainties associated with current emission modeling and impact characterization approaches.

\subsection{Summary of the challenges and limitations of LCA for launch vehicles}
\label{sec:2-lca-summary}
Based on the above state-of-the-art, the main challenges and limitations of current LCA practices for launch vehicles are summarized in Tab.~\ref{tab:limitations-lca}.
\begin{table}[h]
    \caption{\\
    Key challenges and limitations of current LCA practices for launch vehicles.}
    \label{tab:limitations-lca}
    \centering
    \small
    \begin{tabular}{p{0.3\columnwidth} p{0.6\columnwidth}}
        \toprule
        Limitation & Description \\
        \midrule
        \textbf{Space-specific data scarcity} &
        Despite progress (e.g., ESA LCA database, SSSD), high-quality life-cycle inventory data for space systems remain limited, particularly for materials, manufacturing processes, and propellant production. \\
        \midrule
        \textbf{High-altitude emission modelling} &
        Standard climate metrics neglect altitude-dependent effects of rocket emissions. Impacts of species such as aluminum oxide and black carbon, especially on ozone depletion, remain poorly characterized and uncertain. \\
        \midrule
        \textbf{Impact interpretation consistency} &
        Aggregated indicators (e.g., PEF) rely on weighting schemes not tailored to space systems, potentially misrepresenting key categories such as ozone or resource depletion. \\
        \midrule
        \textbf{Late LCA integration} &
        LCA is typically applied post-architecture definition, limiting influence on early design. Early integration within MDAO is needed for meaningful trade-offs. \\
        \bottomrule
    \end{tabular}
\end{table}

This paper addresses the last limitation by integrating predictive LCA within an MDAO framework for the early-stage conceptual design of launch vehicles. The main focus is placed on the methodological development required to enable this integration and, in particular, to support multi-objective optimization involving performance, cost, and environmental criteria. While this study focuses on the framework development and its application, it is designed to enable future extensions, including the integration of the other points discussed (addition of space-specific data, advanced modeling of high-altitude emissions' impacts with accurate characterization factors).

\section{Methodology for including a predictive LCA model within an early-stage MDAO framework for launch vehicles}
\label{sec:3}
\subsection{Description of the MDAO framework}
\label{sec:3-1}
The early-stage design of a launch vehicle is inherently complex due to its strongly coupled multidisciplinary nature, involving disciplines such as propulsion, structure, aerodynamics, weights and sizing, trajectory optimization, and costs. Optimizing them independently leads to a sequential and highly iterative process, as modifications in one discipline propagate to the others. To address this issue, the MDAO framework is employed. 
It enables the simultaneous integration of disciplinary models within a unified optimization architecture, allowing for consistent handling of couplings and constraints while systematically exploring the design space.

More particularly, the general MDAO problem can be formulated as follows \cite{Balesdent2011}:
\begin{align}
    \min \hspace{1cm}& f(\mathbf{z},\mathbf{y}) \label{eq:1}\\
    \text{w.r.t.} \hspace{1cm}& \mathbf{z},\mathbf{y} \notag \\
    \text{s.t.} \hspace{1cm}& \mathbf{g}(\mathbf{z},\mathbf{y}) \leq 0 \label{eq:2}\\
    & \mathbf{h}(\mathbf{z},\mathbf{y}) = 0 \label{eq:3}\\
    & \forall (i,j) \in \{1, \ldots, N\}^2,\; i \neq j, \mathbf{y}_{ij} = \mathbf{c}_{ij}(\mathbf{z}_i, \mathbf{y}_{.i}) \label{eq:4}\\
    & \mathbf{z}_{min} \leq \mathbf{z} \leq \mathbf{z}_{max} \notag
\end{align}
where:
\begin{itemize}
    \item $f(\cdot)$ is the objective function to be optimized (e.g., Gross Lift-Off Weight -- GLOW, payload mass, costs),
    \item $\mathbf{z} \in \mathbb{R}^{n_z}$ is the vector collecting the $n_z$ design variables which may include continuous variables (e.g., propellant masses, stage diameters, mixture ratios) and discrete variables (e.g., propellant mix, material selection) \cite{Baraton2025, Pelamatti2020},
    \item $\mathbf{y}$ is the vector composed of the coupling variables that model the interactions between the disciplines, also being their inputs (e.g., stage dry masses for trajectory, dynamic pressure and axial load factor for structure),
    \item $\mathbf{g}(\cdot) = \begin{bmatrix} g_1(\cdot) & \cdots & g_{n_g}(\cdot) \end{bmatrix}^T$ is the vector of the $n_g$ inequality constraint functions (e.g., target orbit minimum perigee, maximum heat flux, maximum angle of attack),
    \item $\mathbf{h}(\cdot) = \begin{bmatrix} h_1(\cdot) & \cdots & h_{n_h}(\cdot) \end{bmatrix}^T$ is the vector of the $n_h$ equality constraint functions (e.g., target orbit altitude and velocity),
    \item $\mathbf{c}(\cdot)$ is the vector of coupling functions that calculate the coupling variables which come out of a discipline (e.g., specific impulse for propulsion, aerodynamic coefficients for aerodynamics, stage component masses for structure),
\end{itemize}
and $N$ the number of disciplines. 
Note that MDAO problems must be able to handle interdisciplinary coupling satisfaction inside the design process. This requirement is enforced by Eq.\eqref{eq:4} which states that the coupling is said to be satisfied when, for two disciplines $i$ and $j$, the subset of coupling variables $\mathbf{y}_{ij}$ is equal to the output of the subset of coupling functions $\mathbf{c}_{ij}$. Example of such variables include load factors between trajectory and structure, or specific impulse between propulsion and trajectory. 

This general formulation can be further extended to several MDAO problem classes. First, it can be adapted to different MDAO architectures (corresponding to different organization of the solving process of the MDAO problem), including single-level methods such as Multi Discipline Feasible (MDF), Individual Discipline Feasible, and All-At-Once, as well as multi-level methods such as Collaborative Optimization and Bi-Level Integrated System Synthesis. These architectures mainly differ in the way coupling consistency and feasibility are enforced. For more details on the implementation of the different methods and their comparisons, the reader is referred to \cite{Balesdent2011}. 

Another key aspect of multidisciplinary design of launch vehicles is the handling of the trajectory optimization which consists of an optimal control problem. The latter can be solved using numerical methods available in the literature \cite{Betts2010}. Different manners for implementing and solving the trajectory optimization algorithm within the MDAO formulation are possible. Either the optimal control problem is handled inside the trajectory discipline: in that case, the system optimizer does not control the optimal control variables. Or, the optimal control variables, hereafter called $\mathbf{w} \in \mathbb{R}^{n_w}$, are treated at the same level of the system design variables $\mathbf{z}$ by the optimizer. These two approaches are more specifically detailed in \cite{Balesdent2022}. 

Then, once the MDAO architecture is defined, the optimization algorithm to solve the MDAO problem can be either stochastic (e.g., genetic algorithms, evolutionary strategies) or deterministic (e.g., gradient-based approaches) \cite{Martins_EngineeringDesignOpt}. Finally, it can be generalized to a multi-objective optimization problem, where several conflicting criteria are optimized simultaneously. In fact, different optimization strategies can be employed depending on the design goals: classical constrained optimization focuses on finding one or a set of Pareto-optimal solutions while strictly enforcing constraints \cite{Balesdent2022}, whereas quality-diversity optimization approaches aim to explore the design space more broadly, generating a diverse set of feasible solutions that highlight alternative trade-offs \cite{Baraton2025}. Therefore, defining the MDAO problem for launch vehicle design remains a challenging task, and the formulation needs to be adapted to the launch vehicle case study (expendable, reusable, design of the whole launch vehicle or only a stage, etc.).
In the following paragraphs, a formulation of the MDAO problem is defined to illustrate the proposed methodology, however, the authors emphasize that alternative formulations are possible and that the corresponding choice is a topic that remains open for further research.

\begin{figure}[!h]
    \centering
    \includegraphics[width=1\linewidth]{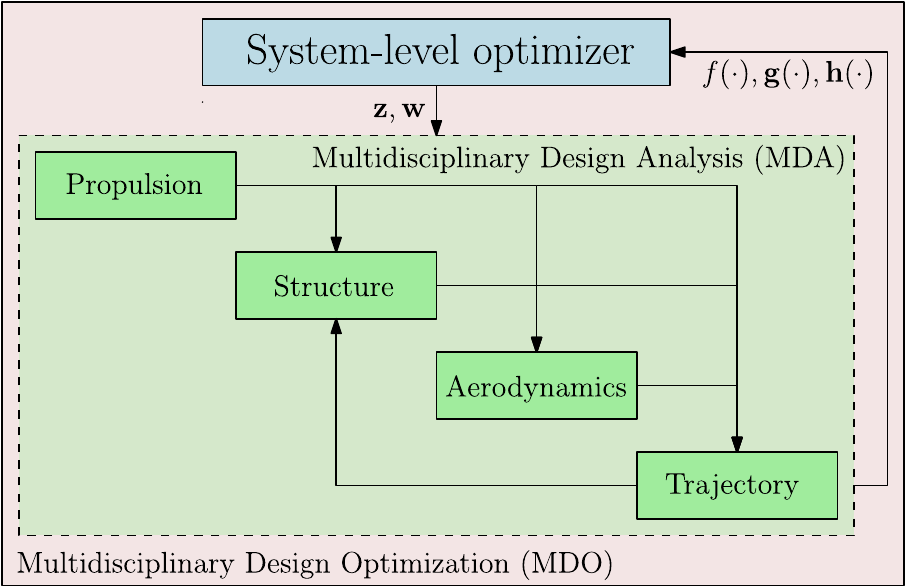}
    \caption{Typical MDF formulation for a launch vehicle design process. In this example, the trajectory optimization is handled at the system level.}
    \label{fig:mdf-legacy}
\end{figure}

In fact, the most used MDAO formulation for launch vehicle design is the MDF method \cite{Balesdent2011}. The latter consists in using a Multidisciplinary Design Analysis (MDA) process to satisfy the coupling consistency between all the different disciplines. 
More particularly, this condition is met when the system of nonlinear equations formed by Eq.~\eqref{eq:4} is satisfied, i.e., when for all disciplines, the coupling variables are consistent across disciplines. 
Thus, all the disciplines are coupled in an analysis module which checks the feasibility of the solution at each iteration. Typically, MDA can be solved using fixed point iteration (e.g., Gauss-Seidel algorithm \cite{Martins_EngineeringDesignOpt}) between all the disciplines. 
Fig.~\ref{fig:mdf-legacy} illustrates a typical example of a MDF formulation for the early design of launch vehicles. 

\begin{figure}[!h]
    \centering
    \includegraphics[width=1\linewidth]{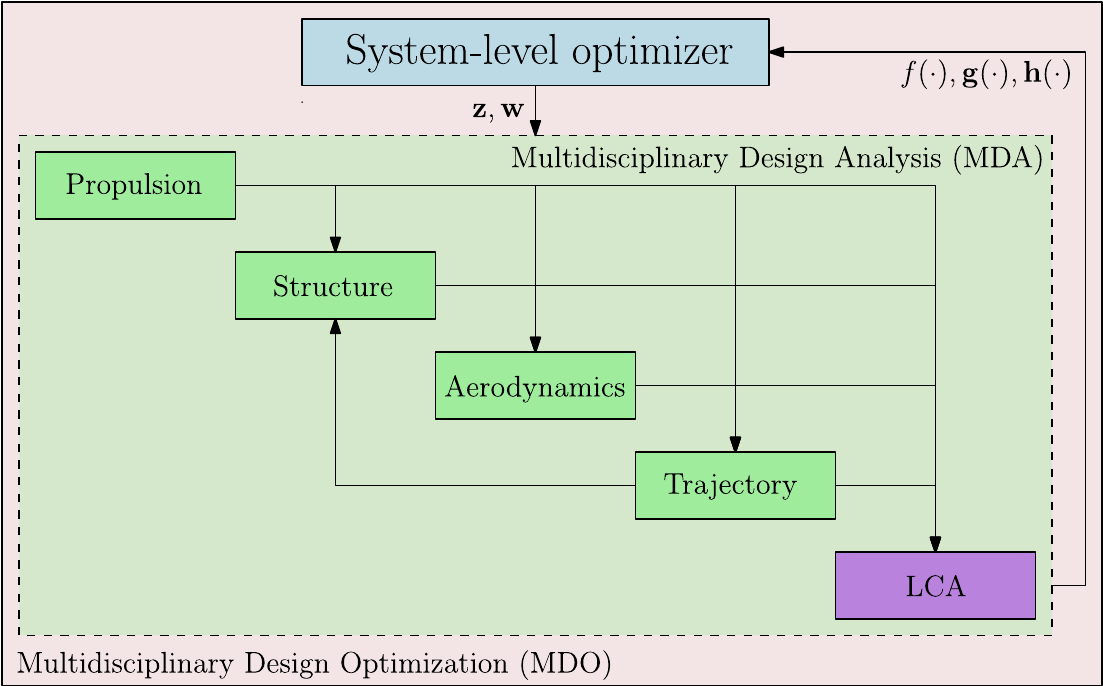}
    \caption{Typical MDF formulation with integration of LCA considerations.}
    \label{fig:mdao_process}
\end{figure}

The main contribution of this paper is the description of a methodology for embedding predictive LCA within an MDAO tool for launch vehicles. It is achieved through the integration of an LCA discipline into the common MDAO framework described above, which is adopted here as one possible architectural choice. 
One can note that alternative MDAO formulations would yield equivalent results, and that this choice is therefore primarily a matter of computational cost and optimization efficiency rather than influencing the outcomes.
Therefore, the MDAO formulation of Fig.~\ref{fig:mdf-legacy} is augmented with the integration of the LCA discipline as illustrated in Fig.~\ref{fig:mdao_process}. 

The versatility of the methodology then offers the possibility of formulating alternative optimization problems to incorporate environmental considerations.
In the context of launch vehicle design, a typical formulation consists in minimizing the GLOW under a set of specification constraints related to the target orbit, as well as physical limitations (e.g., maximum allowable heat flux, angle of attack, or dynamic pressure). Environmental criteria can be incorporated either by including them directly in the objective function, thereby transforming the problem into a multi-objective optimization \cite{Balesdent2025}, or by introducing them as an additional subset of inequality constraints within the function $\mathbf{g}(\cdot)$ \cite{TORMENA2022}. Another possible approach consists in exploring a set of feasible solutions using a quality-diversity optimization approach \cite{Baraton2025}. In this paper, the multi-objective formulation is adopted as an illustrative example of the proposed methodology.

In the next paragraphs, the methodology to integrate the predictive LCA discipline within the MDAO framework is described in details. 

\subsection{Description of the predictive LCA discipline model}
\label{sec:3-2}
Following the ISO 14040:2006 \cite{iso14040}, the discipline must follow the four-step approach mentioned in Sec.~\ref{sec:2-lca-methodology}. 
In this study, in compliance with the ESA guidelines \cite{esa2025lca}, the following functional unit has been selected:
\begin{description}
\centering
    \item[] \it To transport 1 metric ton of payload from the launch site to destination B over 1 year's operations.
\end{description}
Assessing the environmental impacts over one year of operations enables the comparison of different launch vehicle architectures with different payload capacities. For example, one architecture may perform 8 launches carrying 20~t of payload each, whereas another may perform 10 launches carrying 16~t of payload each, while delivering the same total payload over the year (160~t). In this way, the number of launches and the payload capacity of the vehicle are inherently accounted for in the assessment. This functional unit is particularly relevant when performing a comprehensive LCA that includes testing and integration activities as well as reusability capabilities, for which the environmental impacts do not necessarily scale linearly with the number of launches.
Therefore, the results of the LCA are expressed in metric ton of payload delivered to orbit per year (t\textsubscript{payload}.year\textsuperscript{-1}). 

\begin{figure*}[!b]
    \centering
    \includegraphics[width=1\linewidth]{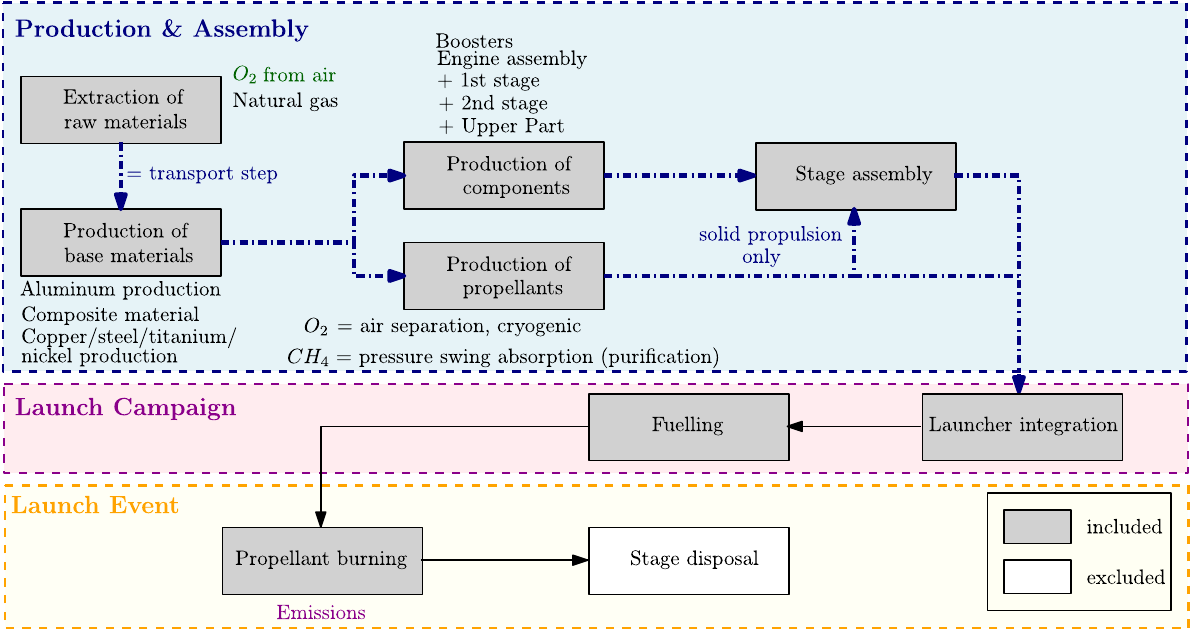}
    \caption{System boundaries for the expendable launch vehicle design, adapted from \cite{esa2025lca}. In gray are represented the phases included in the present study, while in white are represented the phases that should be considered but that are here kept for future works.}
    \label{fig:sys_boundaries}
\end{figure*}

This functional unit implies the definition of system boundaries, defining the scope of the study.
In this paper, the LCA has been restricted to the production and assembly phases of the vehicle until the launch event following the ESA guidelines \cite{esa2025lca}. It is detailed in Fig.~\ref{fig:sys_boundaries}.
More particularly, it includes the production of the stages' components and propellants from raw materials with the transport phase associated.
Additionally, it models the corresponding atmospheric impacts from launch emissions. Nevertheless, environmental impacts during stage disposal are not considered at this stage for simplicity. 
Indeed, the objective of this study is to investigate how predictive LCA can be integrated into the early design phases of launch vehicles, rather than to provide a detailed environmental assessment of an industrial-scale future launch system. Consequently, test and integration activities, as well as recovery and refurbishment processes for reusable launch vehicles, are not included in the present scope. Nevertheless, the proposed methodology is generic and may be extended to account for them by adapting the granularity of the LCA and the system boundaries considered.

\begin{figure}[!h]
    \centering
    \includegraphics[width=1\linewidth]{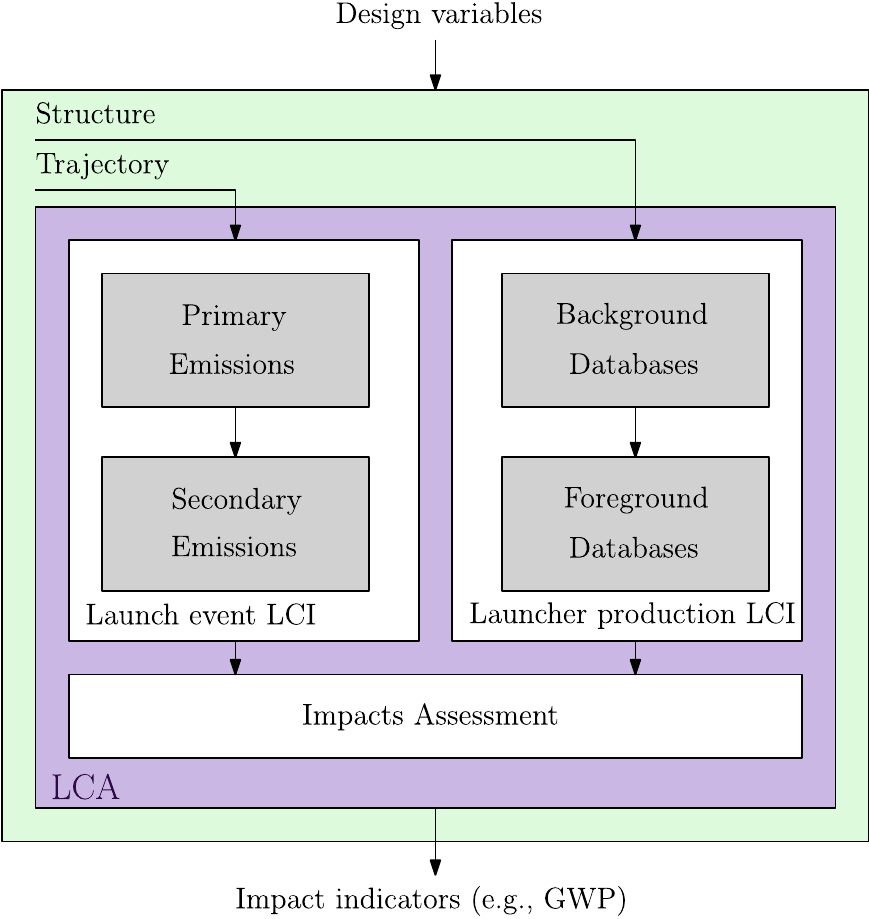}
    \caption{Description of the LCA discipline.}
    \label{fig:LCA-discipline}
\end{figure}

Within the MDAO framework, the challenging task is the inventory analysis since the launcher components' masses vary along the optimization. It is therefore necessary to develop a predictive LCA based on dynamic parameters (here the masses) that are updated at each iterations of the optimization process. Furthermore, it is also necessary to develop a launch vehicle emission model that estimates the emissions and their impact during the ascent trajectory. Therefore, the LCA discipline displayed in Fig.~\ref{fig:mdao_process} is described more precisely in Fig.~\ref{fig:LCA-discipline}.

A three-step approach must be followed to obtain the environmental impacts of the launch vehicle under study: (i) the definition and integration of parametric inventories with respect to the components' masses and the design and coupling variables, (ii) the estimation of the launch vehicle emissions from the trajectory data and the design and coupling variables, (iii) the assessment of the environmental impacts using a characterization model. Note that the first two points can be carried out in parallel as they are not dependent on each other. The methodology is described in the next paragraphs.

\subsubsection{Construction of the parametric inventories}
\label{sec:3-2-inventory}
As mentioned previously, the critical aspect of the methodology defined in this paper is the dynamic inventory analysis involved for each component. This second phase of the LCA methodology relies on the accurate collection of the data related to all the processes that occur throughout the life cycle of the system under study. 
Even if LCI databases exist \cite{ESA_LCA, frischknecht2005ecoinvent, wilson2019thesis} providing aggregated information for unit processes with fully modeled life cycles, identifying the relevant processes and building LCI models by combining them with foreground data derived from in-house primary information remains both resource- and time-consuming. More particularly, because such foreground data can vary with the design being evaluated, for example due to changes in component mass, material, or manufacturing process, the LCI models may need to be appropriately adapted and the corresponding database calls repeated, which also adds a significant computational cost. Consequently, this limits their integration into an MDAO framework.

In this paper, parameterized inventories based on dynamic parameters and supporting fast processing speeds of the impacts assessment are developed. This approach allows the LCA to be continuously updated as a function of the design and coupling variables.
These inventories rely on the definition and integration of foreground and background processes. The background processes, which are not under the direct control of the optimizer (e.g., the production of 1 kg of aluminum), are generated and evaluated prior to their integration into the MDAO framework. In contrast, the foreground processes are defined as algebraic formulas depending on optimization variables from the launcher design phase (e.g., component masses computed within the structural discipline), as well as on the pre-generated background processes.

\begin{figure*}[!h]
    \centering
    \includegraphics[width=1\linewidth]{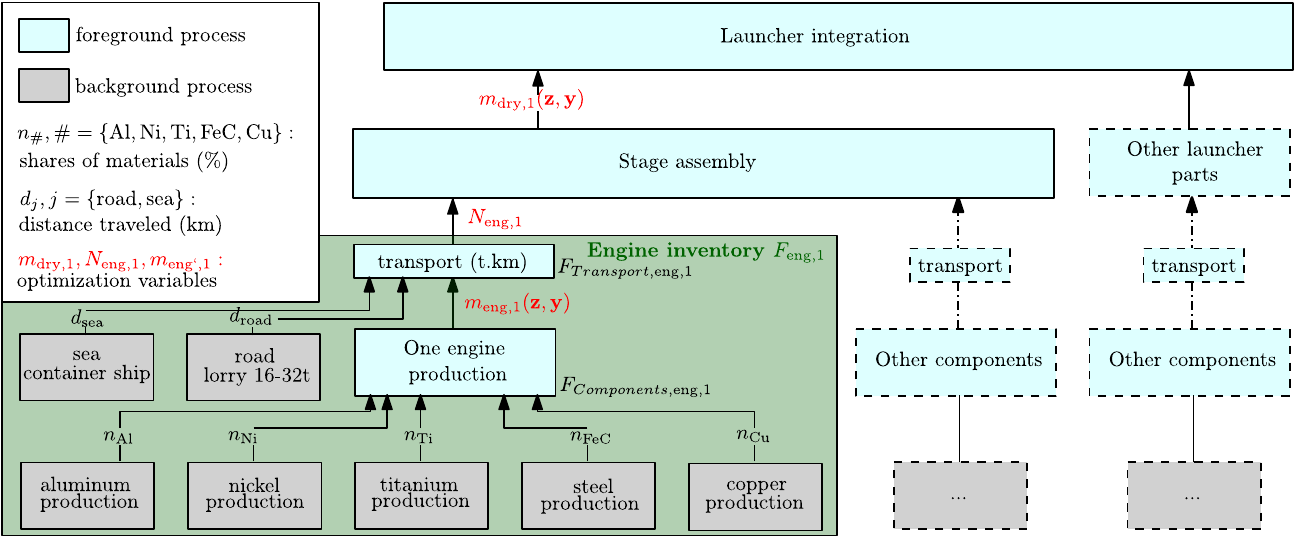}
    \caption{Process tree for the engine production. In gray are represented the background processes (generated via LCI databases), in blue the foreground processes that depend on user's specifications and design or interdisciplinary coupling variables. The latter are depicted in red.}
    \label{fig:dynamic_processtree}
\end{figure*}

Fig.~\ref{fig:dynamic_processtree} provides an overview of the dynamic process tree that is required to fulfill the functional unit taking the example of the first-stage engine production.
As an example, the foreground database, $F_{\mathrm{eng},1}$, denoting the engine production of the first stage as defined in Fig.~\ref{fig:dynamic_processtree} can be expressed as follows:
\begin{equation}
\label{eq:7}
    \begin{aligned}
        F_{\mathrm{eng},1} 
        &= F_{Components,\mathrm{eng},1}(\mathbf{z}, \mathbf{y}) + F_{Transport,\mathrm{eng},1}(\mathbf{z}, \mathbf{y}) \\
        &= m_{\mathrm{eng},1}(\mathbf{z}, \mathbf{y}) \times \bigg[ \big( n_{\mathrm{Al}} \times B_{\mathrm{Al}} + n_{\mathrm{Ni}} \times B_{\mathrm{Ni}} \\
        &\quad + n_{\mathrm{Ti}} \times B_{\mathrm{Ti}} + n_{\mathrm{FeC}} \times B_{\mathrm{FeC}} + n_{\mathrm{Cu}} \times B_{\mathrm{Cu}} \big) \\
        &\quad + \big( d_{\mathrm{road}} \times B_{\mathrm{road}} + d_{\mathrm{sea}} \times B_{\mathrm{sea}} \big) \bigg]
    \end{aligned}
\end{equation} 
where $n_{\#},\# = \{\mathrm{Al, Ni, Ti, FeC, Cu}\}$ for aluminum, nickel, titanium, steel, and copper, and $d_j, j = \{\mathrm{sea}, \mathrm{road}\}$ are user-defined parameters related to the share of materials and the distance of transport (sea or road), respectively. 
One can note that these parameters could also be treated as design or coupling variables, or defined as dependent on them, particularly when discrete design variables reflecting technological or architectural choices are involved (e.g., including the engine cycle feed system as a design variable would modify the material distribution accordingly).
$B_{\#/j}$ are background processes imported from the available LCI databases, and $m_{\mathrm{eng},1}$ is the first stage engine mass computed within the MDAO process. Therefore, the parametric inventory is assigned to dynamic parameters (here the variables of the MDAO process) and is created once and for all.

The parametric engine inventory is subsequently integrated into the parametric inventory of the first stage (see Fig.~\ref{fig:dynamic_processtree}). For instance, when multiple engines are employed, the engine inventory is scaled by the number of engines $N_{\mathrm{eng},1}$, which may be treated as a discrete design variable. More generally, the first-stage inventory is constructed by recombining contributions associated with the dry mass of the first stage $m_{\mathrm{dry},1}$ (coupling variable from the structure discipline), following the same parametric logic. In this inventory structure, the “Stage assembly” process therefore corresponds to the aggregation of the foreground inventories of the components constituting the stage. This approach is consistently applied to all stages and other launcher subsystems until the complete launch vehicle inventory is assembled, with “Launcher integration” similarly representing the aggregation of the stages and other subsystem inventories. Afterward, the propellant contributions are incorporated. For a detailed overview of the launch vehicle inventory defined in this paper and its subsystem decomposition for LCA, the reader is referred to \ref{app:A}. Further refinement of this inventory through a more detailed LCA of the launch vehicle including additional LCI processes (e.g. electricity use, heating) is left for future work.

\subsubsection{Estimation of the launch vehicle emissions}
\begin{figure*}[!b]
    \centering
    \includegraphics[width=0.8\linewidth]{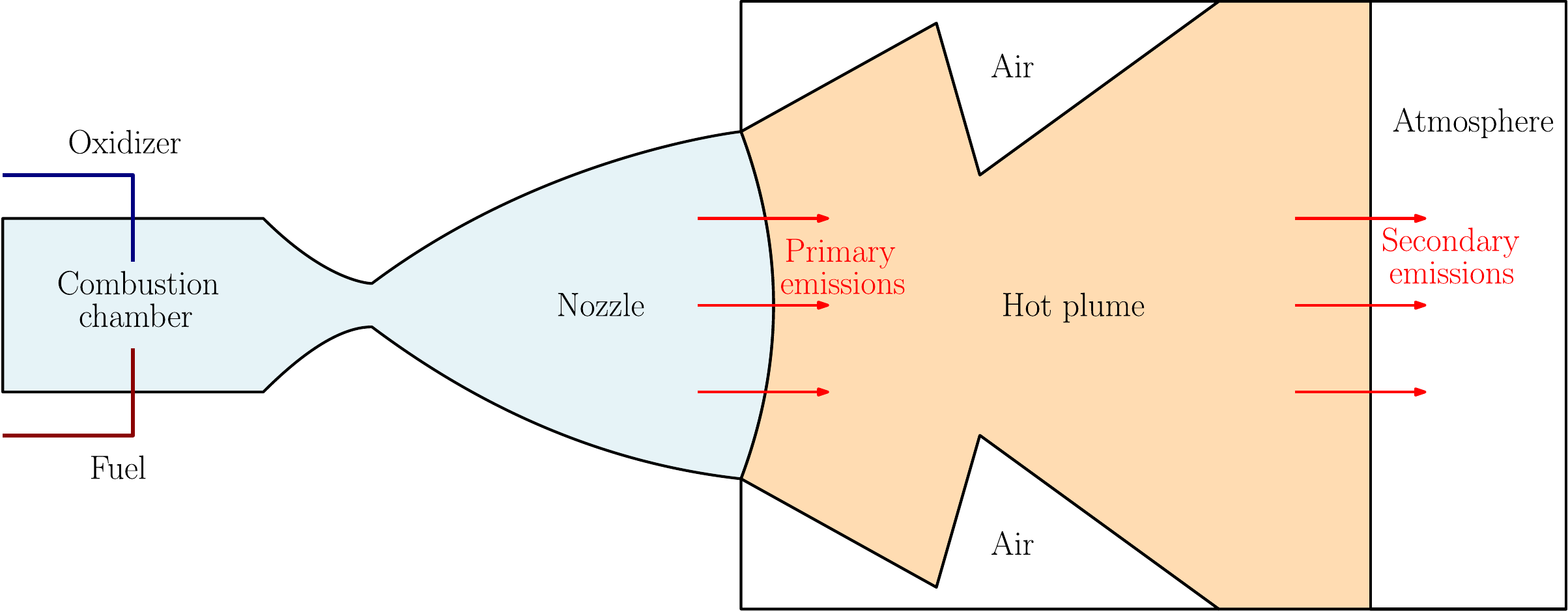}
    \caption{Definition of primary and secondary emissions, adapted from James et al.~\cite{james2021commercial}.}
    \label{fig:launch_emissions}
\end{figure*}

As mentioned in Sec.~\ref{sec:2-lca-sota}, the estimation of the launch emissions is a challenging task and still remains today a critical area for research. 
First, following the literature, let's define the rocket emissions \cite{james2021commercial}. The mass of a species \(i\) emitted by a rocket engine, denoted \(m_{\mathrm{species},i}\), over a trajectory segment of duration \(\delta t = t_f - t_i\) between the initial time \(t_i\) and the final time \(t_f\) (corresponding to an altitude segment between \(h_i\) and \(h_f\)), is given by:  
\begin{equation}
m_{\mathrm{species},i} = \int_{t_i}^{t_f} EI_{\mathrm{species},i}(h(t)) \times \dot{m}(t) \, dt
\end{equation}
where \(EI_{\mathrm{species},i}(\cdot)\) is the emission index associated with the considered species as a function of altitude \(h(\cdot)\), and \(\dot{m}(\cdot)\) is the mass flow rate as a function of time. The emission indices are factors that relate the mass of a species \(i\) emitted by the rocket engine to the mass of propellant consumed. The emission index characterizes the outcome of the chemical reactions involved in the combustion process of the rocket engine for a specific pollutant. In this study, two types of emission indices are considered: the Primary Emission Index (PEI) from the nozzle exit, and the Secondary Emission Index (SEI) resulting from the afterburning reaction between the exhaust plume and the atmosphere. They are generally expressed in $\mathrm{g.kg^{-1}}$. Fig.~\ref{fig:launch_emissions} illustrates them.

Thus, primary emissions are determined by the amount of each species at the nozzle exit. They result from the combustion of gases inside the combustion chamber and the nozzle. These primary emissions do not interact with the atmosphere. They vary depending on the type of propellant and the characteristics of the engine. In this study, the approach adopted to estimate primary emissions is based on the use of the Chemical Equilibrium with Applications (CEA) tool from the NASA \cite{mcbride_gordon_1996_cea}. 
Indeed, in addition to calculating the rocket engine characteristics, NASA CEA can provide the composition of the exhaust gases as well as the properties of the mixture under different assumptions regarding the chemical equilibrium of combustion.
Specifically, it assumes perfect chemical equilibrium in the chamber or frozen composition in nozzle expansions through the Gibbs’s free
energy minimization approach, providing a theoretical maximum performance estimate.
Once calculated, these primary emission indices serve as inputs for the estimation of the secondary emissions.  

Several models exist to calculate secondary emissions. It is important to note here that these models remain highly uncertain today \cite{Brevault2025} and have not been experimentally validated. Currently, three methods are available in a context of early design stages:
\begin{itemize}
    \item Using an analytical model derived from empirical data \cite{desain2014potential, malkin1978environmental}; for example, this is the model by James et al.~\cite{james2021commercial}.
    \item Estimating post-combustion with the local atmosphere using NASA CEA \cite{mcbride_gordon_1996_cea}; this is the method adopted by Fischer et al. \cite{FischerLEATegusphere-2026-3050}.
    \item Building a model based on Computational Fluid Dynamics (CFD) calculations, which directly simulates the dynamics of the launch vehicle plume using a Navier-Stokes solver.
\end{itemize}
However, the latter method is very computationally expensive and therefore not suitable for MDAO studies. One possible solution would be to develop multi-fidelity models, as explained by Balesdent et al.~\cite{Balesdent2025}. In fact, this paper has shown that the model of Fischer et al.~\cite{FischerLEATegusphere-2026-3050} is relatively accurate with respect to the CFD calculations, therefore, it is the model adopted in this study. 
Note that for the estimation of the nitrogen oxids NO$_\mathrm{x}$ and black carbon (BC) secondary emissions, the analytical model of James et al.~\cite{james2021commercial} is employed.
Finally, one can note that the MDA as defined in this paper is generic and compatible with other simulation model of emissions. 

To make it computationally suitable for MDAO simulations which involve repeated function evaluations and therefore repeated calculations of the secondary emissions, surrogate models are created using the Gaussian process regression-based Kriging method \cite{SacksKriging} of the SMT toolbox \cite{saves2024smt}. More particularly, a design of experiments of 500 model calculations has been generated. 
The input variables correspond to the design and coupling variables involving the engine characteristics, as well as the altitude and relative velocity profiles obtained from representative trajectories. 
The output variables are the SEI of the different species involved (e.g., CO$_2$, CO, H$_2$O, H$_2$, BC, NO$_\mathrm{x}$ for this study).

\subsubsection{Impacts assessment}

Recalling Fig.~\ref{fig:LCA-discipline}, once the launch emissions have been estimated and the parameterized inventories constructed, the resulting environmental impacts can be assessed. This is performed in two steps.
First, during initialization, that is the first time the complete inventory of the launch vehicle is created, the impacts of each background processes are calculated. This is the initial calculation that can be time-consuming and that is therefore not suitable for MDAO implementation. The impacts are then substituted into the parametric inventory allowing the creation of algebraic formulas for each impact that involves only the dynamic parameters (recall Eq.~\eqref{eq:7}). Importantly, these algebraic formulas introduce no approximation; they exactly reproduce the results of the full background calculation.
The second step of the impacts assessment, that is the subsequent evaluations of the impacts, can be processed in a very short time and is therefore suitable to be integrated within the framework. An LCIA method must be employed at this stage (e.g., the EF method \cite{Damiani2022}) to generate the environmental impacts to be analyzed, notably the GWP100 or the Ozone Depletion Potential (ODP) scores, to which the impact of launch emissions can be added, as described hereafter.

For the launch emissions, the characterization factors available in the literature only allow their conversion into climate change impact through the GWP100. In fact, there is no characterization factors available with respect to the ODP indicator while launch vehicles do emit ozone-depleting substances, such as gaseous chlorine $\mathrm{Cl_\mathrm{x}}$, nitrogen oxids $\mathrm{NO_\mathrm{x}}$, and hydroxyl radicals $\mathrm{OH}$. Moreover, the available factors remain highly uncertain and do not depend on altitude, which increases uncertainties in the upper atmosphere \cite{DominguezCalabuig2024}. They are described in Tab.~\ref{tab:CF_emissions} for different altitude categories, as reported in \cite{DominguezCalabuig2024} with multiple references here cited; at sea level and for aviation cruise altitude. One can note that there is no data with respect to alumina particles $\mathrm{Al_2O_3}$. 
\begin{table}[h]
    \caption{\\
    GWP100 characterization factors of emissions by species and altitude (kgCO\textsubscript{2}-eq) at Earth's surface and for aviation cruise altitude, adapted from \cite{DominguezCalabuig2024}.}
    \label{tab:CF_emissions}
    \centering
    \small
    \begin{tabular}{c|cl|cl}
        \toprule
        Emission & Ground & Refs. & Aviation & Refs. \\
        \midrule
        CO$_2$ & 1 & \cite{IPCC2014} & 1 & \cite{Lee2021}\\
        CO & 4.0 [1.6--7.6] & \cite{DominguezCalabuig2024, IPCC2014} & 4.0 & \cite{DominguezCalabuig2024}\\
        H$_2$O & 0.0005 & \cite{Sherwood2018} & 0.06 & \cite{Lee2021} \\
        H$_2$ & 12.8 & \cite{Hauglustaine2022} & 12.8 & \cite{DominguezCalabuig2024}\\
        BC & 900 [100--1700] & \cite{DominguezCalabuig2024, IPCC2014} & 1166 & \cite{Lee2021} \\
        NO$_\mathrm{x}$ & 8.5 [7--10] & \cite{DominguezCalabuig2024,Lammel1995} & 114 & \cite{Lee2021}\\
        \bottomrule
    \end{tabular}
\end{table}

The GWP100 impact score of launch emissions is calculated by multiplying the mass of each emitted species within each 1-kilometer altitude segment by its corresponding characterization factor, then summing all the resulting contributions.
This score is then added to the GWP100 obtained for the launcher production and assembly phase inventory. 

Therefore, the methodology described in this section enables the development of an LCA discipline suitable for MDAO integration, allowing environmental considerations to be included at the early design stages of a launch vehicle, alongside other criteria such as overall cost and vehicle performance. This capability is illustrated in the next section through a case study.

\section{Case study: application to a two-stage-to-orbit vehicle}
\label{sec:4}
In this section, the predictive LCA discipline is incorporated into an MDAO framework, as described in Sec.~\ref{sec:3-1} and Fig.~\ref{fig:mdao_process}, which aims to design an expandable launch vehicle. Note that the methodology used to build this discipline is generic and could be adapted to other launcher configuration. The one described in the next paragraph must therefore be considered as an example to illustrate the capabilities of the method for generating informed trade-offs between performance and environmental impact. In fact, the first subsection describes the MDAO framework and displays the result for a single-objective optimization without embedding the LCA discipline. The latter is used as a baseline for the upcoming multi-objective studies, carried out in the next subsections. 

\subsection{Description of the MDAO framework used as test case}
\label{sec:4-1}
The test case considers a TSTO launch vehicle to be designed in order to inject 21 metric tons of payload into a circular and equatorial low Earth orbit of 800 km altitude from the Guiana Space Center. The baseline vehicle is defined as propelled with LOx/LCH\textsubscript{4} for both stages. The first stage is composed of 7 engines and the second stage of a single engine. 
The legacy single-objective MDF formulation of the MDAO problem (considering the coupling satisfaction handle by a MDA) is the following:
\begin{align}
    \mathrm{min} \hspace{1cm}& \mathrm{GLOW}(\mathbf{z}, \mathbf{w}) \label{eq:10}\\
    \mathrm{w.r.t} \hspace{1cm}& \mathbf{z}, \mathbf{w} \notag\\
    \mathrm{s.t.} \hspace{1cm}& \mathbf{g}(\mathbf{z}, \mathbf{w}) \leq 0 \label{eq:12}\\
    & \mathbf{h}(\mathbf{z}, \mathbf{w}) = 0 \label{eq:13}\\
    & \mathbf{z}_{min} \leq \mathbf{z} \leq \mathbf{z}_{max}, \
    \mathbf{w}_{min} \leq \mathbf{w} \leq \mathbf{w}_{max} \notag
\end{align}

The design variables are summarized in Tab.~\ref{tab:design-variables}. For simplicity, the engine characteristics (chamber pressure, expansion ratio) are fixed to standard values of gas generator cycle engines \cite{Sippel2024}. More particularly, we define the chamber pressure at 120~bars for both stages, and we consider an expansion ratio of $\epsilon_1 = 20$ for the first stage and $\epsilon_2 = 120$ for the second stage.
\begin{table}[h]
    \caption{\\
    Design variables of the MDAO framework studied with the bounds considered for optimization.}
    \label{tab:design-variables}
    \centering
    \small
    \begin{tabular}{lccc}
        \toprule
        Variables & Notation & Bounds & Unit\\
        \midrule
        1st stage mixture ratio & $O/F_{1}$ & [1.9, 3.9] & -- \\
        2nd stage mixture ratio & $O/F_{2}$ & [1.9, 3.9] & -- \\
        1st stage propellant mass & $m_{\mathrm{prop},1}$ & [150, 600] & t \\
        2nd stage propellant mass & $m_{\mathrm{prop},2}$ & [5, 200] & t \\
        1st stage mass flow rate & $\dot{m}_{1}$ & [250, 600] & kg/s \\
        2nd stage mass flow rate & $\dot{m}_{2}$ & [5, 400] & kg/s \\
        Diameter (equal for both stages) & $D$ & [3.5, 4] & m \\
        Trajectory commands & $\mathbf{w}_{}$ & -- & -- \\
        \bottomrule
    \end{tabular}
\end{table}
The total number of design variables is 13 (3 for each stage, 1 common to both stages, and 6 for the trajectory optimization). More details about the trajectory optimization discipline and the associated design variables are provided in the following paragraphs. Then, seven constraints are present in the optimization problem: the altitude of the apogee of the transfer orbit must be equal to 800 km; the altitude of the perigee of the transfer orbit must be greater than 250 km; the angle of attack of the vehicle during the ascent phase must never exceed 15 deg (with specific restrictions for the atmospheric flight); the axial load factor must be lower than 5 g; the dynamic pressure must be lower than 50 kPa; the heat flux must be lower than 100 $\mathrm{W.m^{-2}}$; and finally, the remaining propellant mass before the circularization burn must be sufficient to circularize the final orbit 800$\times$800 km.

The legacy disciplinary models are adapted from the MDAO framework FELIN (Framework for Evolutive Launcher optImizatioN) available online\footnote{L. Brevault et al. FELIN: Framework for Evolutive Launcher Optimization. \url{https://github.com/l-brevault/FELIN}} and based on the opensource framework OpenMDAO \cite{gray2019openmdao}. In fact, the propulsion discipline is modified with NASA CEA \cite{mcbride_gordon_1996_cea} to compute propulsion-specific data. Furthermore, the structure discipline is augmented with more detailed mass models \cite{castellini2013multidisciplinary}, especially for the second stage of the launch vehicle. 
\begin{description}
   \item[Propulsion discipline:] The disciplinary model is developed using NASA CEA \cite{mcbride_gordon_1996_cea}. 
   Based on the selected propellant type (e.g., LOx/LCH\textsubscript{4}), fixed engine parameters such as the chamber pressure and the nozzle expansion ratio, as well as the MDAO design variables, namely the oxidizer-to-fuel ratio \(O/F\) and the mass flow rate \(\dot{m}\), the model estimates for each stage, the vacuum specific impulse \(I_{sp,vac}\) and the nozzle exit area \(A_e\).
   In fact, NASA CEA delivers theoretical performance predictions, making it well-suited for conceptual and preliminary design stages. The resulting specific impulse and nozzle exit area values act as coupling parameters and are propagated to both the structure and the trajectory disciplines.

   \item[Structure discipline:] This discipline is adapted from the geometry and weights estimation models as defined by Castellini \cite{castellini2013multidisciplinary} and enables the computation of the dry masses of the launch vehicle stages. These models rely on weight estimation relationships derived from statistical analysis of a database of existing launchers.
   Depending on its architecture, the launch vehicle can be decomposed into solid propellant boosters, lower stage and intermediate stages. 
   Without loss of generality, for the particular case study of this paper, only two stages are considered.
   The first stage accounts for the engine assembly mass (engine, thrust vector control, pressurization system), the thrust frame, the fuel and oxidizer tanks (as well as the intertank if necessary), the thermal protection systems associated, and the interstage. Then, the second stage considers the avionics and the electrical power system modules, the fairing and payload adapter, and all the other components mentioned for the first stage, except the interstage.  
   Therefore, several inputs are required to compute the stage masses. For each stage \(i=\{1,2\}\), these include the number of engines \(N_{\mathrm{eng},i}\) (here fixed); the MDAO design variables, namely the propellant mass \(m_{\mathrm{prop},i}\), the oxidizer-to-fuel ratio \(O/F_i\), the mass flow rate \(\dot{m}_i\), and the stage diameter \(D_i\); as well as structural loads for the launch vehicle, such as the maximum axial load factor \(n_{x,\max}\) and the maximum dynamic pressure \(Q_{\mathrm{dyn},\max}\), which are either specifications (and handled as constraints in the optimization problem) or coupling variables provided by the trajectory discipline.
   The outputs of the discipline are composed of the dry masses of each stage $m_{\mathrm{dry},i}$, necessary for the trajectory discipline. Furthermore, it includes the detailed mass decomposition between all the components as inputs for the predictive LCA discipline (see \ref{app:A}).

   \item[Aerodynamics discipline:] The aerodynamics discipline focuses on computing aerodynamic coefficients, such as drag and lift, which are essential for evaluating aerodynamic loads during the launcher’s atmospheric flight. Estimating the aerodynamic performance across all flight regimes (subsonic, transonic, supersonic, and hypersonic) is a complex and demanding task. In this approach, a previously calculated drag coefficient table is directly integrated in the discipline. The resulting drag force scales according to geometrical parameters such as the stage diameter. The latter is computed offline using low fidelity aerodynamics code, similar to MissileDATCOM \cite{MissileDATCOM2014}, for a representative two-stage-to-orbit expendable launcher. Therefore, the resulting drag coefficient table, provided as function of Mach number $M$ and angle of attack $\alpha$, hereafter defined as $C_x(M,\alpha)$, is directly supplied to the trajectory discipline. This eliminates the need for an iterative feedback loop between the aerodynamics and trajectory models. Such an approach is generally adequate for early-phase design studies. 

   \item[Trajectory discipline:] The trajectory discipline consists of integrating the equations of motion over time for the launch ascent flight. 
   The latter is assumed to happen in the equatorial plane for a non-rotating Earth and is modeled in 2D polar coordinates. The jettisoning of the first stage and the payload fairing for the TSTO vehicle is considered and constrained according to threshold values of dynamic pressure and heat flux. The payload is injected into a circular orbit using a Hohmann transfer ascent. During the trajectory, aerodynamic loads, drag, and lift are computed as functions of velocity and altitude, based on the U.S. Standard Atmosphere models \cite{US1976}. The considered equations of motion are: 
   \begin{align}
        \dot{r} &= V \, \sin\gamma, \\
        \dot{\lambda} &= \frac{V \cos\gamma}{r}, \\
        \dot{V} &= \frac{-D + T \cos(\theta-\gamma)}{m} - g \sin\gamma, \\
        \dot{\gamma} &= \left(\frac{V}{r} - \frac{g}{V}\right) \cos\gamma + \frac{T \sin(\theta-\gamma)}{m V}, \\
        \dot{m} &= - q,
   \end{align}
   where $r$ is the position vector from Earth's center to a fixed point in the vehicle (e.g., the center of gravity), $\lambda$ is the longitude, $\theta$ is the pitch angle, $\gamma$ is the flight path angle, $\alpha$ is the angle of attack, $m$ is the mass of the vehicle at any given time, $q$ is the propellant mass flow rate, $g$ is the acceleration due to gravity, $D$ is the drag force, $T$ is the thrust force, and $V$ is the velocity of the launch vehicle.
   The integration is done using a 5th-order Runge-Kutta method with a dedicated structure on top to accounts for different discrete events and their appearance order such as fairing jettison, changes in the control law based on flight conditions, and other mission-specific transitions. A direct single-shooting method is used to define the optimal control law, represented by a parameterized profile of the pitch angle \cite{castellini2013multidisciplinary}. 
   The launch is divided into four flight phases: lift-off, pitch-over maneuver, gravity turn, and bi-linear tangent law. A discontinuity in the control law is permitted between the gravity turn and bi-linear tangent phases, as aerodynamic forces become negligible at higher altitudes. Each phase is parameterized by a specific set of control variables. During the lift-off phase controlled by a duration $\Delta t_{LO}$, the pitch angle is set equal to 90 deg imposing a vertical flight. Then, a linear pitch-over maneuver is carried out, controlled by three parameters: the duration of the phase $\Delta t_{PO}$, the maximum variation of the pitch angle $\Delta \theta_{PO}$, and the duration of the exponential decay $\Delta t_{decay}$. More particularly, this maneuver is composed of two phases: first the pitch angle is linearly increased up via the law $\theta = \gamma - \Delta \theta_{PO} \cdot t / \Delta t_{PO}$, followed by an exponential decay phase to match the gravity turn conditions ($\theta = \gamma$). This is enforced via the law $\theta(t) = \gamma - \Delta \theta_{PO} \cdot \exp{ \big(-t / \Delta t_{decay} \big)}$. After three time constants $\Delta t_{decay}$, these conditions are assumed reached and a gravity turn flight is performed with an angle of attack equal to zero (no control parameters are enforced during this phase). Once a lower threshold on the dynamic pressure is reached, this maneuver is followed by an exo-atmospheric phase that is parameterized by a bi-linear tangent law: $\tan(\theta) = \frac{a^{\xi} \, \tan\theta_i + \big(\tan\theta_f - a^{\xi} \, \tan\theta_i\big) t}{a^{\xi} + \big(1 - a^{\xi}\big) t}$, where $a$ is an arbitrary constant, $\xi$ is a control parameter defining the shape of the bi-linear tangent law, and $\theta_i$ and $\theta_f$ are the pitch angles at the beginning and at the end of the phase, respectively, also considered as control variables. Finally, the thrust is throttled in the flight of the first stage in order to respect the maximal axial load factor constraint.
   For further details on the guidance for launch vehicle in MDAO context, one can consult \cite{castellini2013multidisciplinary}.
   This discipline requires as inputs the design variables, the coupling variables obtained via the other disciplines (e.g., specific impulse, stage dry masses) as well as the trajectory command parameters for the pitch angle control law defined above, composing the vector $\mathbf{w} \in \mathbb{R}^{n_w}$. More particularly, these parameters defined as design variables in the MDAO problem (recall Sec.~\ref{sec:3-1}), are summarized in Tab.~\ref{tab:traj-variables}. For simplicity, the duration of the exponential decay, $\Delta t_{decay}$, is set equal to 10 s, and the arbitrary constant of the bi-linear tangent law, $a$, is set equal to 100.
   It outputs the trajectory profiles of the launch ascent flight, which are necessary for the predictive LCA discipline to estimate the launch emissions.
   \begin{table}[h]
        \caption{\\
        Control variables for trajectory optimization, $\mathbf{w}$, with their respective bounds. In the current MDAO framework, they are considered as design variables.}
        \label{tab:traj-variables}
        \centering
        \small
        \begin{tabular}{lcc}
            \toprule
            Variables & Bounds & Unit\\
            \midrule
            $\Delta t_{LO}$ & [5, 20] & s \\
            $\Delta t_{PO}$ & [5, 20] & s \\
            $\Delta \theta_{PO}$ & [1, 10] & deg \\
            $\xi$ & [-1, 1] & -- \\
            $\theta_i$ & [0, 50] & deg \\
            $\theta_f$ & [-20, 20] & deg \\
            \bottomrule
        \end{tabular}
   \end{table}
\end{description}

\begin{figure*}[!b]
    \centering
    \includegraphics[width=0.9\linewidth]{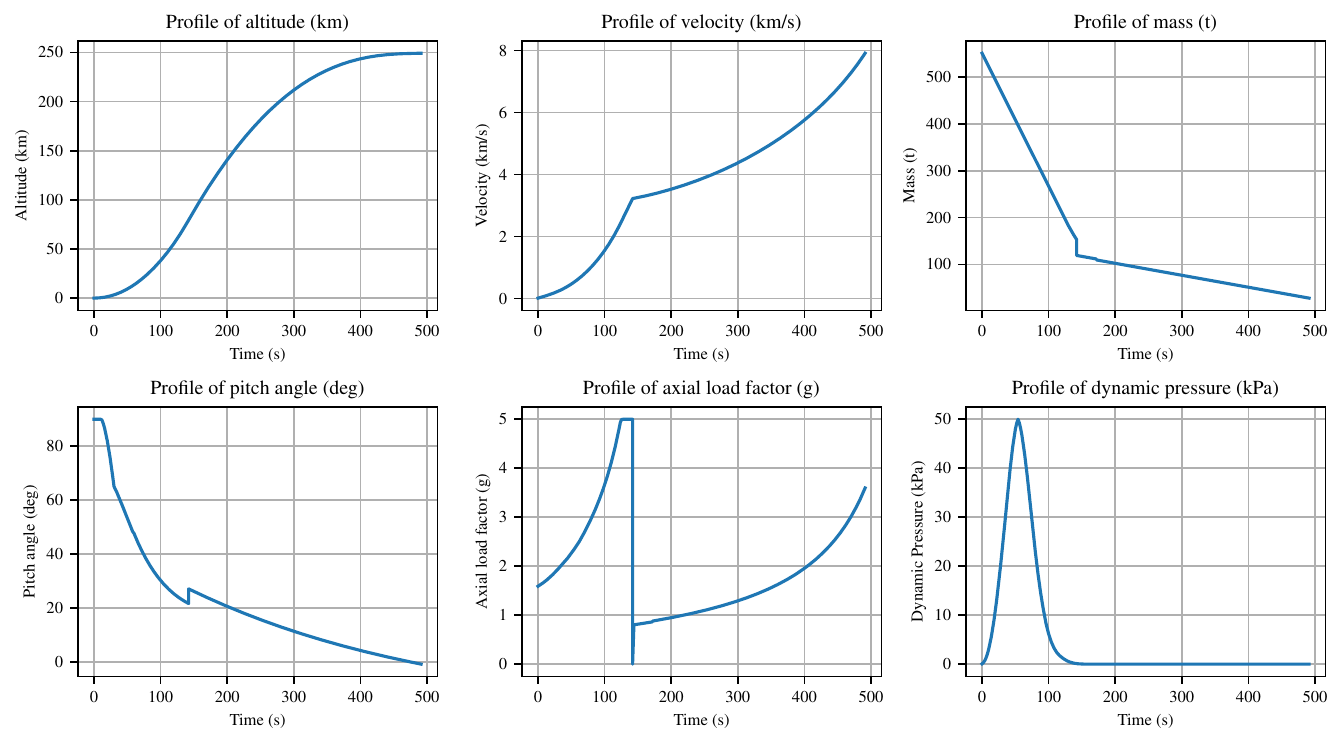}
    \caption{Ascent trajectory for the baseline (min GLOW) TSTO launch vehicle (propelled phase).}
    \label{fig:baseline-traj}
\end{figure*}

The Covariance Matrix Adaptation - Evolution Strategy (CMA-ES) algorithm \cite{hansen2001cmaes} is used to solve the optimization problem. 
The resulting optimal design variables and the main features of the TSTO launch vehicle when only optimizing the GLOW and without considering any environmental impact criterion, are given in Tab.~\ref{tab:baseline-res}. 
\begin{table}[!h]
    \caption{Optimal design variables and main characteristics of the baseline (minimum GLOW) TSTO launch vehicle.}
    \label{tab:baseline-res}
    \centering
    \small
    \begin{tabular}{lc|lc}
        \toprule
        \multicolumn{2}{l}{\textit{Objective:}} & GLOW (t) & 550.85 \\
        \midrule
        \multicolumn{2}{l|}{\textit{Design variables:}} & \multicolumn{2}{l}{\textit{Coupling variables:}} \\
        \midrule
        $O/F_{1}$ & 2.88 & $I_{sp,vac,1}$ (s) & 339.47 \\
        $O/F_{2}$ & 2.96 & $I_{sp,vac,2}$ (s) & 362.80 \\
        $m_{\mathrm{prop},1}$ (t) & 397.82 & $A_{e,1}$ (m$^2$) & 1.27 \\
        $m_{\mathrm{prop},2}$ (t) & 89.30 & $A_{e,2}$ (m$^2$) & 4.79 \\
        $\dot{m}_{1}$ (kg\,s$^{-1}$) & 406.37 & $m_{\mathrm{dry},1}$ (kg) & 33977.66 \\
        $\dot{m}_{2}$ (kg\,s$^{-1}$) & 255.50 & $m_{\mathrm{dry},2}$ (kg) & 6750.27 \\
        $D$ (m) & 4.00 & Stage 1 length (m) & 42.01 \\
         &  & Stage 2 length (m) & 11.01 \\
        \bottomrule
    \end{tabular}
\end{table}
The optimal trajectory obtained for this configuration is displayed in Fig.~\ref{fig:baseline-traj} with the evolutions of altitude, relative velocity, mass, pitch angle, dynamic pressure and axial load factor. Note that the optimal trajectory injects the payload at the perigee of the transfer orbit. 

On this baseline configuration, the environmental impact indicators can then be computed using the predictive LCA discipline described in the previous section. 
The opensource framework Brightway2 \cite{mutel2017brightway} is employed to handle LCI databases and LCIA methods. It enables the use of the lca\_algebraic library \cite{jolivet2021lca} for constructing the parametric inventory of the launch vehicle under study.
The latter was constructed following the methodology of Sec.~\ref{sec:3-2-inventory} as described in details in \ref{app:A}. 
Regarding the definition of the functional unit, since an expendable launch vehicle is considered and R\&D, testing, and integration activities are neglected, the annual formulation of the functional unit (t\textsubscript{payload}.year\textsuperscript{-1}) is equivalent to a normalization by the payload capacity (t\textsubscript{payload}\textsuperscript{-1}), which is here 21~t. This is because both the environmental impacts and the delivered payload over the year scale linearly with the number of launches. Moreover destination B corresponds to the transfer orbit, as mentioned above.

In this study, the LCIA method selected is the EF method involving 16 indicators \cite{Damiani2022}. With respect to space activities, in addition to climate change and ozone depletion already mentioned, the most critical impacts are the resources depletion using the Abiotic Depletion Potential (ADP) for fossil fuels as well as minerals and metals, the land use, and the water depletion \cite{sala2017weighting, Verkammen2023eucass}. Moreover, as mentioned in Sec.~\ref{sec:2-lca-mdao}, the European Commission also introduced a normalization and weighting procedure yielding the PEF single score, expressed in points (pts). More specifically, the latter is obtained as follows:
\begin{equation}
    \begin{aligned}
    \label{eq:15}
    \mathrm{PEF}(\mathbf{z}, \mathbf{w}) &= \sum_{k \in \mathcal{E}} \mathrm{PEF}_k(\mathbf{z}, \mathbf{w}) 
    = \sum_{k \in \mathcal{E}} \left( \frac{k(\mathbf{z}, \mathbf{w})}{\mathrm{NF}_k} \times \mathrm{WF}_k \right), \\
    &\text{with } \mathcal{E} = \{\mathrm{GWP}100, \mathrm{ODP}, \mathrm{ADP\ (fossil)}, \\
    & \hspace{1.5cm} \mathrm{ADP\ (metals)}, \mathrm{Land\ use}, \mathrm{Water\ use}\}
    \end{aligned}
\end{equation}
where $\mathrm{NF}$ corresponds to the normalization factor of the EF impact category, and $\mathrm{WF}$ its weighting factor. The values for each impact category are available in ~\cite{Damiani2022}.
The results are summarized in Tab.~\ref{tab:baseline-GWP}. 
\begin{table}[!h]
    \caption{\\
    Environmental impacts of the baseline (minimum GLOW) TSTO launch vehicle following the EF method, expressed in t\textsubscript{payload}\textsuperscript{-1}. The last column shows the contribution of the impact to the aggregated single score.}
    \label{tab:baseline-GWP}
    \centering
    \small
    \begin{tabular}{lcc}
        \toprule
        Environmental indicators & Values & \% PEF \\
        \midrule
        GWP100 (kgCO$_2$-eq) & 276284.43 & 18.82 \\
        Water use (m$^3$world-eq) & 982658.53 & 18.25\\
        ADP fossil fuels (MJ) & 2450407.96 & 7.84 \\
        ADP minerals/metals (kgSb-eq) & 16.75 & 54.84 \\
        Land use (pts) & 1189904.91 & 0.18\\
        ODP (kgCFC-11-eq) & 0.01 & 0.07\\
        \midrule
        PEF single score (pts) & 3983.28 & 100 \\
        \bottomrule
    \end{tabular}
\end{table}

Once applying the normalization and weighting procedure, it can be observed that this single score is mainly driven by the impact of resources depletion of metal and mineral components (more than 50\% of the total score). It is then followed almost equally by the impacts on water use and climate change (18-19\%, respectively). Finally, the impact on resources depletion of fossil fuels is around 8\% of the single score, whereas the impacts on land use and ozone depletion are marginal. For the sake of conciseness, GWP100 will hereafter be referred to as GWP.

To assess the order of magnitude of the results obtained with the proposed framework, a comparison with available LCA studies of launch vehicles was performed. However, it is important to note that direct quantitative comparisons remain challenging because the reported results depend strongly on the functional unit, system boundaries, vehicle architecture, life-cycle stages considered, and assumptions used for the different life-cycle processes. Moreover, the level of detail provided in the available literature is generally insufficient to fully reconstruct the underlying life-cycle inventories and modeling assumptions. The comparison should therefore be regarded as an order-of-magnitude consistency check rather than a formal validation of the model. 

First, Tab.~\ref{tab:gwp_score_comparison} compares the GWP value of the baseline (minimum GLOW) solution of this study with values of LCA studies available in the literature for different launch vehicle types. 
\begin{table}[!h]
    \centering
    \caption{\\
    Comparison of GWP values, expressed in t\textsubscript{payload}\textsuperscript{-1}, with previous launch vehicle LCA studies.}
    \label{tab:gwp_score_comparison}
    \small
    \begin{tabular}{cccc}
        \hline
        Study & Vehicle type & Propellant & 
        GWP (tCO$_2$-eq) \\
        \hline
        Baseline
        & TSTO
        & LOx/LCH$_4$ 
        & \makecell[c]{
            276\\
            \footnotesize (Emissions: 89)
        } \\
        
        \cite{Harris2019}
        & Falcon 9 
        & LOx/RP-1 
        & $\sim$2470 \\
        
        \cite{Harris2019}
        & Falcon Heavy 
        & LOx/RP-1 
        & $\sim$880 \\

        \cite{DominguezCalabuig2022}
        & TSTO VTVL
        & LOx/LCH$_4$
        & \makecell[c]{
            $\sim$1400$^{*}$\\
            \footnotesize (Emissions: $\sim$191$^{*}$)
        } \\
        \hline
    \end{tabular}
    
    \vspace{2mm}
    \begin{minipage}{0.95\textwidth}
        \footnotesize
        $^{*}$ Normalized from the reported 28,000 tCO$_2$-eq for a 20 t payload.
    \end{minipage}
\end{table}
More particularly, Harris and Landis \cite{Harris2019} assessed the GWP of Falcon 9 and Falcon Heavy rockets. It is worth noting that the contribution of launch emissions reported in their assessment is particularly low (approximately 1\%), which may reflect differences in the assumptions and modeling approach used to account for launch emissions.
Then, Dominguez Calabuig et al.~\cite{DominguezCalabuig2022} assessed reusable launch vehicle types with different fuels (hydrogen, ammonia, kerosene, methane and propane) and first-stage recovery options. Among them, we selected for comparison the LOx/LCH$_4$ Vertical Take-Off Vertical Landing (VTVL) configuration since it is the closest architecture to the case study of this paper. For consistency with the present study, the GWP values reported in these studies were converted, when necessary, to a common basis of metric tons of CO$_2$-eq per metric ton of payload delivered to orbit. The GWP obtained for this study is of the same order of magnitude as the values reported in these previous assessments, although its lower value may indicate a slight underestimation of the overall GWP, which could also result from differences in system boundaries and modeling assumptions.

Additional LCA studies of Ariane launch vehicles were also considered, including the assessments of Ariane 64 by Gallice et al.~\cite{gallice2018environmental}, respectively, as well as the more recent assessment of Ariane 62 by Pellegrino \cite{Pellegrino2025}. Although these studies do not provide sufficiently detailed absolute GWP values for a direct numerical comparison, they provide information on the relative contribution of the different life-cycle stages. However, these results should be interpreted with caution, as the Ariane configurations include solid propellants, whose production and use contribute significantly to the overall GWP and therefore limits the direct comparability with the LOx/LCH$_4$ configuration considered in this study. More particularly, the LCA analysis of the case study is divided into 4 stages: the production of components, the synthesis of propellants, the transport associated, and the launch emissions. The comparison is depicted in Tab.~\ref{tab:gwp_contributions_comparison}. 
\begin{table}[!h]
    \centering
    \caption{\\
    Comparison of main GWP contributions with previous launcher LCA studies.}
    \label{tab:gwp_contributions_comparison}
    \small
    \begin{tabular}{cccc}
        \hline
        Study & Vehicle type & Propellant & \% of total GWP \\
        \hline
        Baseline
        & TSTO
        & LOx/LCH$_4$ 
        & \makecell[l]{
        Components: 35\%; \\
        Propellants: 33\%; \\
        Emissions: 32\%; \\
        Transport: 0\%  
        } \\
        
        \cite{gallice2018environmental}
        & Ariane 64 
        & \makecell[c]{
        LOx/LH$_2$ \\
        + solid }
        & \makecell[l]{
        Propellants: $\sim$45\%; \\
        Components: $\sim$30\%; \\
        Launch$^{*}$: $\sim$25\%; \\
        Transport: $\sim$1\% 
        } \\
        
        \cite{Pellegrino2025}
        & Ariane 62 
        & \makecell[c]{
        LOx/LH$_2$ \\
        + solid }
        & \makecell[l]{
        Components: $\sim$20\%; \\
        Launch$^{\dagger}$: $\sim$55\%; \\
        Transport: $\sim$15\%; \\
        Testing: $\sim$10\% 
        } \\
        \hline
    \end{tabular}

    \vspace{2mm}
    \begin{minipage}{0.95\linewidth}
        \footnotesize
        $^{*}$ Launch campaign included but launch emissions neglected.\\
        $^{\dagger}$ Launch includes propellant production and first stage propellant use (via primary emissions only, black carbon excluded).
    \end{minipage}
\end{table}
Gallice et al.~\cite{gallice2018environmental} reported contributions of approximately 45\% from propellant production and 30\% from production and assembly of the components, while Pellegrino \cite{Pellegrino2025} reported a dominant contribution from the launch event (approximately 55\%), followed by manufacturing (20\%), transport (15\%), and testing (10\%). However, in addition with the different propellant types, Gallice et al.~\cite{gallice2018environmental} neglected the launch emissions, while Pellegrino \cite{Pellegrino2025} considered the primary emissions of the first stage only and accounted for the propellant production within the "Launch" life-cycle stage. Despite these differences in vehicle architecture, system boundaries, and modeling assumptions, the relatively even distribution observed in the present study between component production (35\%), propellant production (33\%), and launch emissions (32\%) is broadly consistent with the range of contributions reported in the literature. It is also worth noting that Pellegrino \cite{Pellegrino2025} reported a similar distribution of the impact categories contributing to the PEF single score, with ADP for metals and minerals representing approximately 40\% of the total score and GWP approximately 25\%, which is consistent with the dominance of ADP for metals and minerals observed in the present study.

\subsection{Multi-objective optimizations}
For the upcoming studies, the predictive LCA discipline is integrated in the optimization loop and the MDAO problem is solved, for the same test case, as a multi-objective optimization considering both the GLOW and an user-defined environmental indicator among the EF method impacts categories or directly the PEF single score as objective functions to be minimized. 
More particularly, the optimization problem becomes as follows:
\begin{align}
    \mathrm{min} \hspace{1cm} & \Big\{ 
    \mathrm{GLOW}(\mathbf{z}, \mathbf{w}), f_{\text{LCIA}}(\mathbf{z}, \mathbf{w}) \Big\}, \label{eq:MDAO_obj} \\
    \mathrm{w.r.t.} \hspace{1cm} & \mathbf{z}, \mathbf{w} \notag \\
    \mathrm{s.t.} \hspace{1cm} & \mathbf{g}(\mathbf{z}, \mathbf{w}) \leq 0 \label{eq:MDAO_ineq} \\
    & \mathbf{h}(\mathbf{z}, \mathbf{w}) = 0 \label{eq:MDAO_eq} \\
    & \mathbf{z}_{\min} \leq \mathbf{z} \leq \mathbf{z}_{\max}, \quad 
      \mathbf{w}_{\min} \leq \mathbf{w} \leq \mathbf{w}_{\max} \notag
\end{align}
where $f_{\text{LCIA}}(\cdot, \cdot)$ is equal to a user-defined indicator among the EF method impacts categories (e.g., $\mathrm{GWP}100(\cdot, \cdot)$) or to the single score $\mathrm{PEF}(\cdot, \cdot)$ as defined in Eq.\eqref{eq:15}.
The design variables are still defined according to Tab.~\ref{tab:design-variables} in Sec.~\ref{sec:4-1} and the constraints are unchanged. 

In this paper, to illustrate the range of analyses that can be performed, we focus on three indicators: the climate change via the GWP100, the water depletion expressed as cubic meters of water use related to the local water scarcity (m$^3$world-eq), and the aggregated PEF single-score indicator defined in Eq.\eqref{eq:15}. The methodology could be extended to a higher number of environmental impact indicators into a three or more objectives, however, in order to analyze preliminary trade-offs, a bi-objective problem is considered. 
The multi-objective optimization was conducted using a hybrid approach. It is first performed using the Comma-Selection Multi-Objective Covariance Matrix Adaptation - Evolution Strategy (COMO-CMA-ES) algorithm \cite{COMO-CMA-ES}, configured with 10 CMA-ES kernels and a population size of 24 individuals, with a maximum of 30,000 iterations. Then, a set of representative non-dominated solutions were selected via a K-medoids clustering method \cite{maranzana1963, park2009} to seed the genetic algorithm NSGA-II \cite{NGSA2}. The latter was run with a population size of 40 individuals for 15,000 generations. No explicit numerical convergence criterion was imposed; instead, a sufficiently large iteration budget was prescribed, and the resulting Pareto front was assessed based on its stabilization and coverage. This sequential strategy was adopted to combine the exploration capabilities of COMO-CMA-ES with the subsequent refinement and densification of the Pareto front using NSGA-II.

The optimization was performed on an HPC cluster using 24 CPU cores in parallel. A 30,000-iteration COMO-CMA-ES run required approximately 50 h of wall-clock time, while the subsequent NSGA-II optimization required approximately 25 h. The total computational time was therefore approximately 75 h for the two-stage optimization. The computational workflow was not specifically optimized for execution time, and further optimization of the parallelization and computational workflow could reduce this computational cost.

It is important to emphasize that the trends observed hereafter are inherently dependent on the modeling assumptions, the definition of the case study (e.g., launcher architecture and propellant choice), and the scope and level of detail of the LCA (e.g., type of propellant production, transport, and manufacturing processes). Consequently, the results presented in this paper should be interpreted as indicative of relative trends within this specific framework.

\subsubsection{Minimizing GLOW and GWP}
\begin{figure}[!h]
    \centering
    \includegraphics[width=1\linewidth]{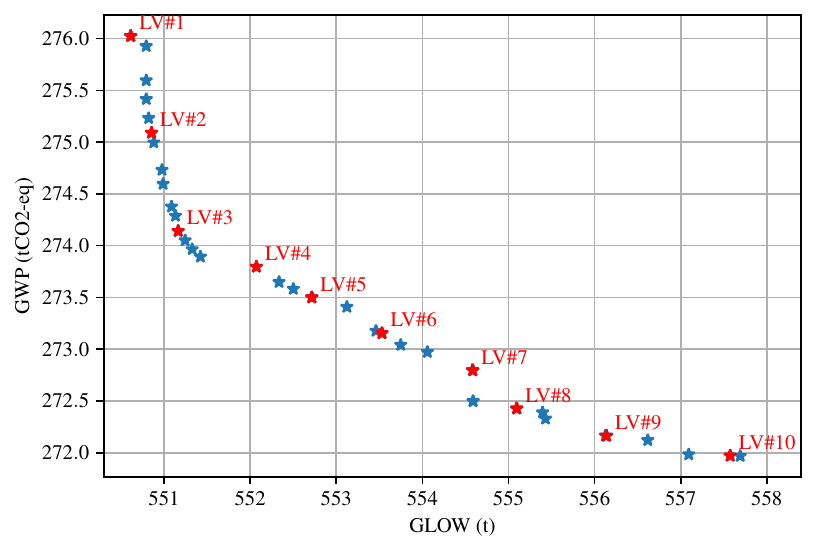}
    \caption{Pareto front between GWP and GLOW. In blue are represented the full set of solutions. Among them, 10 solutions are highlighted in red for the upcoming analyses.}
    \label{fig:GWP_pareto}
\end{figure}

In this subsection, the EF method impact category indicator used is the GWP for assessing the impacts on the climate change. The latter is expressed in tCO\textsubscript{2}-eq as mentioned in Sec.~\ref{sec:3-2}. From the MDAO framework developed in this study, it includes the impacts from the launch vehicle production and assembly, as well as from the flight emissions (see Fig.~\ref{fig:sys_boundaries}). The obtained Pareto front is displayed in Fig.~\ref{fig:GWP_pareto}.

To analyze the set of feasible solutions, the 10 solutions highlighted in red in the figure are studied more deeply. The latter are obtained through a K-medoids clustering method \cite{maranzana1963, park2009} and therefore are considered as the most representative subset of solutions of the Pareto front obtained. First of all, the two extreme solutions in the Pareto front are compared: LV\#1, i.e., the one having the minimum GLOW (thus maximum GWP), and LV\#10, i.e., the one having the minimum GWP (thus maximum GLOW). Tab.~\ref{tab:optima-comparison} shows the differences among the design variables and therefore the disciplinary or coupling variables.
\begin{table}[!h]
    \caption{Comparison of the minimum GLOW and minimum GWP solutions. 
    Relative difference $\Delta$ (\%) corresponds to the minimum GWP solution 
    relative to the minimum GLOW solution.}
    \label{tab:optima-comparison}
    \centering
    \small
    \begin{tabular}{lccc}
        \toprule
         & Min GLOW & Min GWP & \multirow{2}{*}{$\Delta$ (\%)} \\
         & LV\#1 & LV\#10 & \\
        \midrule
        \multicolumn{4}{l}{\textit{Objectives:}} \\
        \midrule
        GLOW (t) & 550.6 & 557.6 & 1.3 \\
        GWP (tCO$_2$-eq) & 276.0 & 272.0 & $-$1.5 \\
        \midrule
        \multicolumn{4}{l}{\textit{Design variables:}} \\
        \midrule
        $O/F_{1}$ & 2.9 & 3.2 & 10.8 \\
        $O/F_{2}$ & 3.0 & 3.0 & 2.8 \\
        $m_{\mathrm{prop},1}$ (t) & 397.1 & 415.1 & 4.5 \\
        $m_{\mathrm{prop},2}$ (t) & 89.9 & 78.8 & $-$12.3 \\
        $\dot{m}_{1}$ (kg/s) & 405.8 & 409.7 & 1.0 \\
        $\dot{m}_{2}$ (kg/s) & 255.6 & 220.3 & $-$13.8 \\
        $D$ (m) & 4.0 & 4.0 & 0.0 \\
        \midrule
        \multicolumn{4}{l}{\textit{Coupling variables:}} \\
        \midrule
        $I_{sp,vac,1}$ (s) & 339.5 & 337.7 & $-$0.5 \\
        $I_{sp,vac,2}$ (s) & 362.8 & 362.6 & $-$0.1 \\
        $A_{e,1}$ (m$^{2}$) & 1.3 & 1.3 & 0.5 \\
        $A_{e,2}$ (m$^{2}$) & 4.8 & 4.1 & $-$13.9 \\
        $m_{\mathrm{CH_4},1}$ (kg) & 102423.7 & 99151.2 & $-$3.2 \\
        $m_{\mathrm{CH_4},2}$ (kg) & 22702.2 & 19500.5 & $-$14.1 \\
        $m_{\mathrm{LOx},1}$ (kg) & 294646.2 & 315954.0 & 7.2 \\
        $m_{\mathrm{LOx},2}$ (kg) & 67140.1 & 59294.0 & $-$11.7 \\
        $m_{\mathrm{dry},1}$ (kg) & 33930.3 & 34576.0 & 1.9 \\
        $m_{\mathrm{dry},2}$ (kg) & 6770.3 & 6098.8 & $-$9.9 \\
        1st stage length (m) & 41.9 & 42.8 & 2.1 \\
        2nd stage length (m) & 11.1 & 9.9 & $-$10.4 \\
        \bottomrule
    \end{tabular}
\end{table}

The comparison between the minimum GLOW and minimum GWP solutions highlights how environmental criteria can reshape launcher architecture. In the minimum GWP configuration, the first stage operates at a higher oxidizer-to-fuel ratio, resulting in reduced methane consumption and a greater share of liquid oxygen. This shift is driven by the lower GWP associated with LOx production and its comparatively lower emissions' impacts during flight. Consequently, first-stage propellant mass increases slightly, and the stage becomes marginally longer and heavier, with corresponding increases in oxidizer tank and related structural masses. 
In contrast, the second stage shows a clear reduction in mass flow rate and in propellant mass along with decreases in tank, engine, and subsystem masses. This leads to a shorter and lighter upper stage with reduced thrust and a smaller nozzle exit area. 
Overall, the minimum GWP solution shifts propulsive effort and structural mass slightly toward the first stage while lightening the second stage and reducing fuel consumption. However, this configuration requires approximately 7 metric tons of additional total propellant, which directly explains the observed increase in GLOW. 
It should be noted that for both optimal configurations (and, more generally, for all optimized solutions), the stage diameter converges to 4 meters, corresponding to the upper bound of the design variable (see Tab.~\ref{tab:design-variables}). This behavior indicates that increasing the diameter systematically reduces structural mass and thus improves the objective function. In fact, in the present formulation of the MDAO framework, no physical or operational constraint counterbalances the structural mass benefit associated with a larger diameter. Consequently, the optimizer naturally drives this variable to its upper bound.

In addition, it is interesting to look at the GWP distribution with respect to the life-cycle stages. Tab.~\ref{tab:GWP-details} shows the results in absolute values and as percentage of the total GWP by life-cycle stage.
\begin{table}[h]
    \caption{\\
    Comparison of the GWP distribution for the minimum GLOW, the median, and the minimum GWP solutions. Values are shown both as absolute GWP and as percentages of total GWP by life cycle stage. Relative difference $\Delta$ (\%) corresponds to the minimum GWP solution relative to the minimum GLOW solution. 
    }
    \label{tab:GWP-details}
    \centering
    \small
    \begin{tabular}{lcccc}
        \toprule
          & Min GLOW & Median & Min GWP & \multirow{2}{*}{$\Delta$ (\%)} \\
          & LV\#1 & LV\#5 & LV\#10 &  \\
        \midrule
        GWP (tCO\textsubscript{2}-eq) & & & & \\
        \midrule
        Emissions & 88.8 & 87.1 & 85.7 & $-$3.5 \\
        Propellants  & 91.4 & 90.7 & 90.7 & $-$0.8 \\
        Components  & 95.6 & 95.4 & 95.3 & $-$0.3 \\
        Transport & 0.2 & 0.2 & 0.2 & 0.0 \\
        \midrule
         \% of total GWP & & & & \\
        \midrule
        Emissions & 32.2 & 31.9 & 31.5 & -- \\
        Propellants & 33.1 & 33.2 & 33.3 & -- \\
        Components & 34.6 & 34.9 & 35.1 & -- \\
        Transport & 0.1 & 0.1 & 0.1 & -- \\
        \bottomrule
    \end{tabular}
\end{table}

Except the transport which plays a marginal role, all the other categories almost contribute equally to the total GWP. But interestingly, 75\% of the GWP reduction of the minimum GWP solution compared to the minimum GLOW solution comes from the launch emissions. 
Increasing the mixture ratio, and therefore reducing methane consumption enables a significant decrease in the climate change impact of launch emissions. 
Conversely, producing more oxidizer in the minimum GWP solution does not increase the impact from the propellant production life-cycle stage, since oxidizer production has a lower climate change impact than methane production. 
Similarly, the climate change impact associated with the production of the launch vehicle components is slightly reduced. This reduction is not solely related to the small decrease in total dry mass (0.06\% relative difference), but also to the redistribution of dry mass between the stages and their different LCI models. In particular, the minimum GWP solution involves a lighter second stage, whose LCI model includes components such as avionics and the electrical power system (see \ref{app:A}), which contribute significantly to the GWP impact. This observation highlights the importance of carefully defining the LCI models within the LCA discipline and ensuring that component masses that are not directly driven by the design variables, such as avionics and EPS in this study, remain physically consistent with the considered vehicle configuration.
In terms of GWP distribution, the lower part of Tab.~\ref{tab:GWP-details} shows that minimizing GWP reduces the share attributed to launch emissions and therefore increases the relative contribution of the other stages. The contribution from propellant production slightly increases, while the component production stage becomes the dominant contributor.

\begin{figure*}[!b]
    \centering
\begin{adjustbox}{minipage=\linewidth,scale=1,center}
    \includegraphics[height=7.8cm]{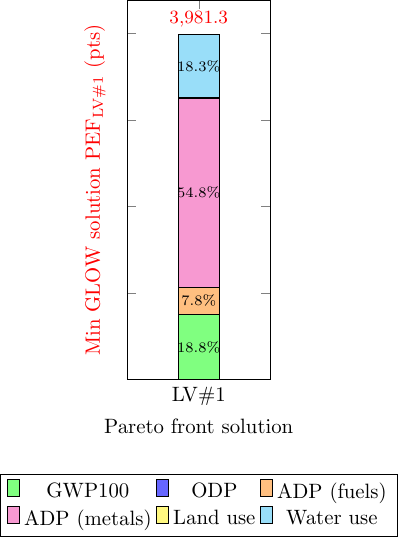}%
    \hfill
    \includegraphics[height=8cm]{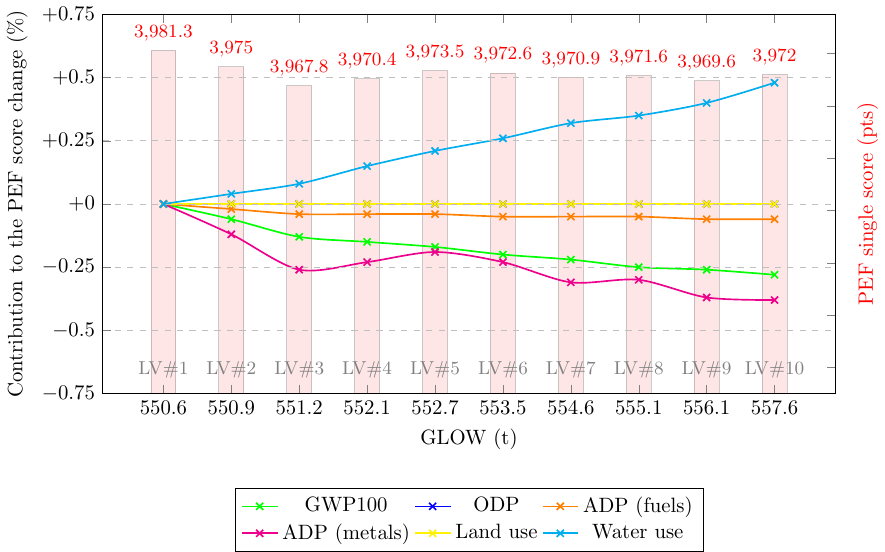}%
\end{adjustbox}
\caption{Left: Contribution of each impact category to the PEF single score for the min GLOW solution ($\mathrm{PEF}_{\mathrm{LV}\#1}$). Right: Evolution of the PEF single score within the multi-objective optimization of GLOW and GWP (background) and contribution of each impact category to the PEF score change with respect to the min GLOW solution (foreground).}
    \label{fig:GWP-PEF}
\end{figure*}

Finally, it is also important to recall that the climate change is not the only impact to consider when conducting an LCA. 
Consequently, when looking for a tradeoff between GLOW and GWP, it is relevant to look at the consequences on the other impact categories and especially on the single score. 
Therefore, Fig.~\ref{fig:GWP-PEF} displays the detailed evolution of the PEF single score for the set of feasible solutions from the Pareto front obtained when carrying out the multi-objective optimization of the GLOW and the GWP. 
More specifically, the contributions of each impact category to the change in the PEF score of the considered solution are highlighted with corresponding colors. These contributions indicate which impact category is responsible for the increase or decrease of the overall PEF score compared to the minimum GLOW solution. They are computed as follows for each impact category: 
\begin{equation}
    \begin{aligned}
    \Delta \mathrm{PEF}_{k, \mathrm{LV}\#i} \big|_{\mathrm{LV}\#1} 
    = \frac{\mathrm{PEF}_{k, \mathrm{LV}\#i} - \mathrm{PEF}_{k, \mathrm{LV}\#1}}{\mathrm{PEF}_{\mathrm{LV}\#i}} \times 100, \\
    k \in \mathcal{E},\\
    i = 2, \dots, 10
    \end{aligned}
\end{equation}
where $k \in \mathcal{E}$ represents a user-defined indicator among the EF method impacts categories (see Eq.\eqref{eq:15}) and LV\#$i$ denotes one of the ten feasible solutions highlighted in Fig.~\ref{fig:GWP_pareto}. On the left side of the figure and similarly to Tab.~\ref{tab:baseline-GWP}, the distribution of the impact category contributions to the PEF score is shown for the minimum GLOW solution, taken as the reference. This provides a baseline for interpreting the contributions to the PEF score changes, which are displayed for each impact category on the right side of the figure.

From this graph, it can be observed that minimizing the GWP with respect to the GLOW is not directly correlated with the PEF single score, or more generally, with other environmental impacts. In fact, this optimization significantly increases the impact on water depletion, which is primarily due to the higher production of liquid oxygen as an oxidizer, a process that consumes substantial amounts of water. More specifically, for most of the solutions, the contribution of the water use impact to the PEF score change is greater in absolute value than that of the climate change impact, meaning that the reduction in GWP does not compensate for the increase in water use. As a result, in terms of the PEF score, these solutions are often environmentally degraded.
Another interesting behavior is the evolution of the impact on resource depletion of metals and minerals, expressed via ADP. This impact is directly related to the dry mass of the launch vehicle and its staging strategy. However, minimizing GLOW and GWP under the current modeling assumptions, is primarily driven by variations in the oxidizer-to-fuel ratio of the first stage, as well as the propellant mass and mass flow rate of the second stage, which translates into reduced launch emissions. These design choices affect the launch vehicle structural dimensions and staging but not necessarily reduces the dry mass, which explains why minimizing GWP is not directly correlated with minimizing ADP for metals and minerals. It should be noted that these conclusions are highly dependent on the mass models adopted in this study, as well as the LCA background processes selected within the available databases (see \ref{app:A}).
In contrast, the ADP for fossil fuels does show a correlation, as it is directly linked to methane production. 
Overall, this analysis highlights the antagonistic behaviors among the different environmental impacts that can arise when performing eco-design for launch vehicles. It also underscores the crucial importance of the single-score definition, since results are highly sensitive to the weighting and normalization procedures used.

\subsubsection{Minimizing GLOW and water use}
\begin{figure}[!h]
    \centering
    \includegraphics[width=1\linewidth]{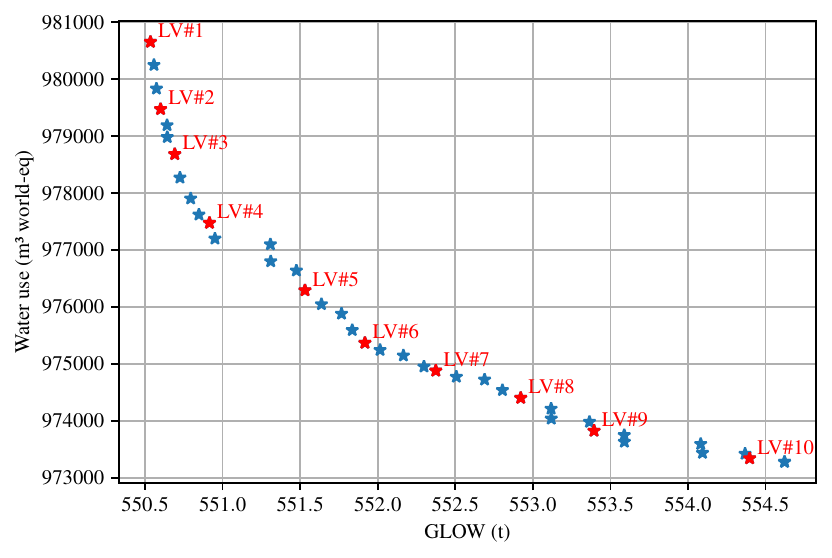}
    \caption{Pareto front between water use and GLOW. In blue are represented the full set of solutions. Among them, 10 solutions are highlighted in red for the upcoming analyses.}
    \label{fig:water_pareto}
\end{figure}

To better understand the antagonistic behaviors, this subsection considers water use as the impact category to be minimized alongside the GLOW, still via the EF LCIA method.  The obtained Pareto front is displayed in Fig.~\ref{fig:water_pareto}. 
Water use is expressed in m$^3$ world-eq, which corresponds to the amount of water consumed in cubic meters weighted by local water scarcity. This unit accounts not only for the volume of water used, but also for its location, giving greater importance to water use in water-scarce regions. From this figure, it can be observed that the range of GLOW values explored during the minimization of water use is smaller than that obtained during the minimization of GWP, suggesting a stronger correlation between these two objectives. 
Tab.~\ref{tab:glow-water-comparison} shows the differences among the design variables and the coupling variables for the extreme solutions like in the analysis of the previous subsection.
\begin{table}[!h]
    \caption{Comparison of the minimum GLOW and minimum water use solutions. 
    Relative difference $\Delta$ (\%) corresponds to the minimum water use solution 
    relative to the minimum GLOW solution.}
    \label{tab:glow-water-comparison}
    \centering
    \small
    \begin{tabular}{lccc}
        \toprule
         & Min GLOW & Min water use & \multirow{2}{*}{$\Delta$ (\%)} \\
         & LV\#1 & LV\#10 & \\
        \midrule
        \multicolumn{4}{l}{\textit{Objectives:}} \\
        \midrule
        GLOW (t) & 550.5 & 554.4 & 0.7 \\
        Water use (m$^3$world-eq) & 980650.9 & 973342.9 & $-$0.7 \\
        \midrule
        \multicolumn{4}{l}{\textit{Design variables:}} \\
        \midrule
        $O/F_{1}$ & 2.9 & 2.6 & $-$8.0 \\
        $O/F_{2}$ & 3.0 & 2.9 & $-$1.7 \\
        $m_{\mathrm{prop},1}$ (t) & 400.0 & 405.1 & 1.3 \\
        $m_{\mathrm{prop},2}$ (t) & 87.0 & 85.7 & $-$1.6 \\
        $\dot{m}_{1}$ (kg/s) & 403.6 & 407.2 & 0.9 \\
        $\dot{m}_{2}$ (kg/s) & 246.3 & 236.8 & $-$3.9 \\
        $D$ (m) & 4.0 & 4.0 & 0.0 \\
        \midrule
        \multicolumn{4}{l}{\textit{Coupling variables:}} \\
        \midrule
        $I_{sp,vac,1}$ (s) & 339.5 & 338.2 & $-$0.4 \\
        $I_{sp,vac,2}$ (s) & 362.8 & 362.8 & 0.0 \\
        $A_{e,1}$ (m$^{2}$) & 1.3 & 1.3 & 0.5 \\
        $A_{e,2}$ (m$^{2}$) & 4.6 & 4.4 & $-$3.9 \\
        $m_{\mathrm{dry},1}$ (kg) & 33927.3 & 34169.1 & 0.7 \\
        $m_{\mathrm{dry},2}$ (kg) & 6599.1 & 6481.8 & $-$1.8 \\
        1st stage length (m) & 42.2 & 43.5 & 3.0 \\
        2nd stage length (m) & 10.8 & 10.7 & $-$0.9 \\
        GWP (tCO$_2$-eq) & 275.3 & 281.1 & 2.1 \\
        \bottomrule
    \end{tabular}
\end{table}

Indeed, the relative difference between the objectives of the minimum water use solution and the minimum GLOW solution is small. The main driver for minimizing water use is the oxidizer-to-fuel ratio of the first stage. As expected from the previous analysis, in contrast to GWP minimization, reducing water use requires a lower $O/F$ ratio to decrease liquid oxygen consumption, which consequently increases the share of liquid methane. For reference, the GWP value for both solutions is also reported; it increases along the optimization, confirming its antagonistic behavior with respect to water use. Only minor adjustments are then needed for propellant masses and mass flow rates. Additionally, a small but noticeable shift of structural mass toward the first stage is again observed.
Overall, this indicates that the design space for minimizing water use is very close to that of the minimum GLOW solution, and the trade-off between launcher performance and water use is relatively weak. 

\begin{figure*}[!h]
    \centering
\begin{adjustbox}{minipage=\linewidth,scale=1,center}
    \includegraphics[height=7.8cm]{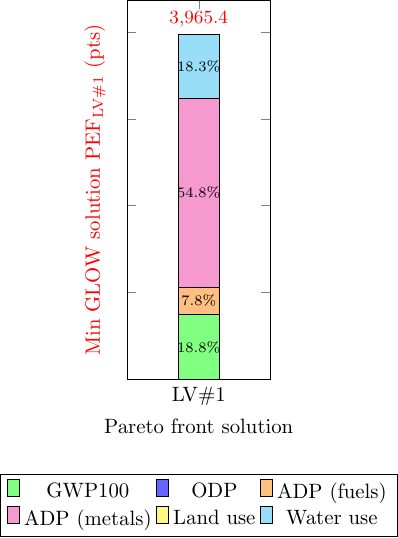}%
    \hfill
    \includegraphics[height=8cm]{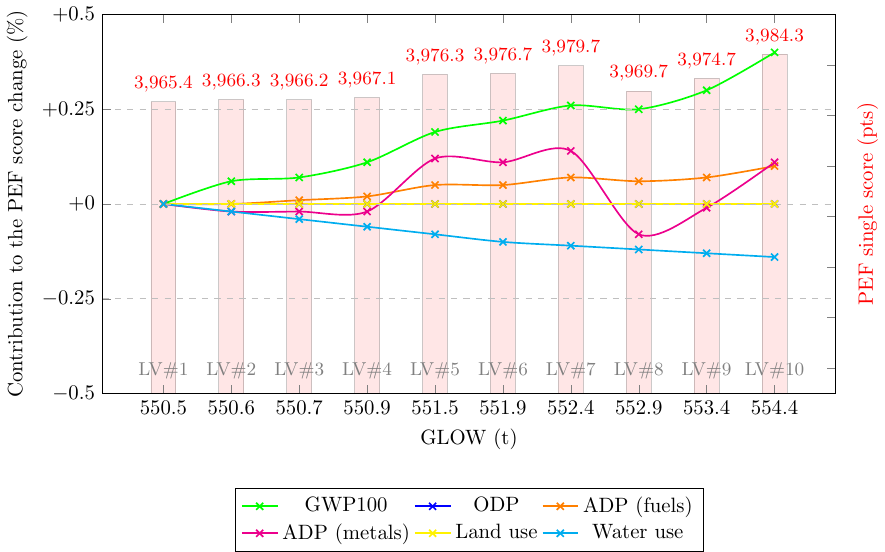}%
\end{adjustbox}
\caption{Left: Contribution of each impact category to the PEF single score for the min GLOW solution ($\mathrm{PEF}_{\mathrm{LV}\#1}$). Right: Evolution of the PEF single score within the multi-objective optimization of GLOW and water use (background) and contribution of each impact category to the PEF score change with respect to the min GLOW solution (foreground).}
    \label{fig:water-PEF}
\end{figure*}

However, as mentioned previously, such a behavior is clearly antagonistic compared to the other environmental impact categories as well as the PEF single-score. This is observed in Fig.~\ref{fig:water-PEF} 
that, similarly to Fig.~\ref{fig:GWP-PEF}, displays the evolution of the single score within the multi-objective optimization and the contributions of each impact category to the PEF score change with respect to the minimum GLOW solution. 

In fact, the figure shows that while minimizing the water use environmental impact, most of the other environmental impacts increase. This is particularly the case for the GWP indicator for climate change, as well as ADP for fossil fuels depletion, due to the higher share of liquid methane, as mentioned previously. 
Regarding the depletion of metal and mineral resources, the ADP indicator does not exhibit a monotonic evolution with respect to GLOW or water use.
In this study, ADP for metals and minerals is closely linked to the vehicle dry mass and its distribution between stages.
During water use minimization, the observed reduction in the second stage propellant mass and mass flow rate indicates a transfer of the second stage towards the first stage. This shift leads to a relative lightening of second-stage components, notably the engine which contains critical materials (e.g., nickel, copper, titanium). Moreover, as mentioned in the previous subsection, the second stage as defined in this case study embeds the electrical power system and avionics, which are major contributors to ADP due to their high material intensity and reliance on critical elements (e.g., lithium-ion batteries). No scaling constraints are imposed on these subsystems with respect to the overall vehicle. 
As a result, for the initial candidate architectures, ADP increases with GLOW due to rising dry mass, but beyond a certain point (LV\#4), mass redistribution from the second to the first stage reduces the contribution of material-intensive subsystems, leading to a variation in ADP. 
More generally, due to its large contribution, these variations strongly influence the PEF single score. 

\subsubsection{Minimizing GLOW and PEF single score}
\begin{figure}[!h]
    \centering
    \includegraphics[width=1\linewidth]{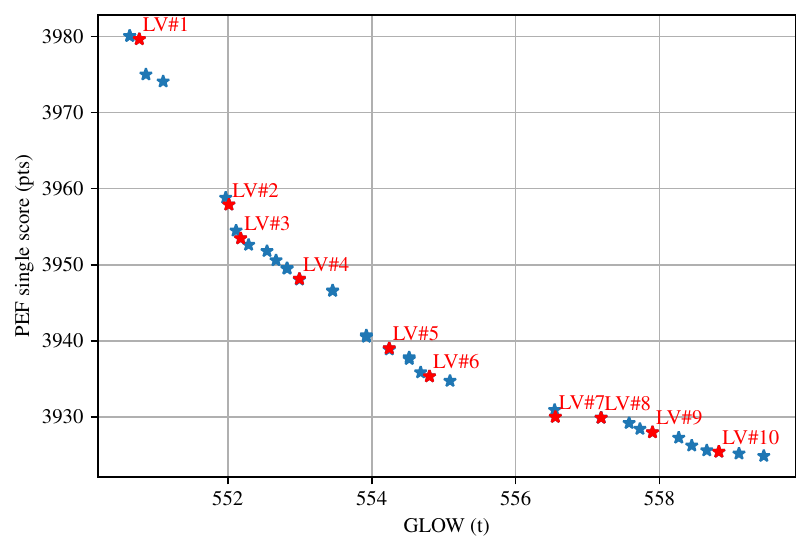}
    \caption{Pareto front between PEF and GLOW. In blue are represented the full set of solutions. Among them, 10 solutions are highlighted in red for the upcoming analyses.}
    \label{fig:SS_pareto}
\end{figure}

Finally, in the following analysis, the optimizer is set to minimize both the PEF single score and GLOW, in order to identify the potential trade-offs, based on the score’s current definition. The obtained Pareto front is displayed in Fig.~\ref{fig:SS_pareto} while Tab.~\ref{tab:glow-ss-comparison} shows the differences among the design variables and the coupling variables for the extreme solutions. 
\begin{table}[!h]
    \caption{Comparison of the minimum GLOW and minimum PEF solutions. 
    Relative difference $\Delta$ (\%) corresponds to the minimum PEF solution 
    relative to the minimum GLOW solution.}
    \label{tab:glow-ss-comparison}
    \centering
    \small
    \begin{tabular}{lccc}
        \toprule
         & Min GLOW & Min PEF & \multirow{2}{*}{$\Delta$ (\%)} \\
         & LV\#1 & LV\#10 &  \\
        \midrule
        \multicolumn{4}{l}{\textit{Objectives:}} \\
        \midrule
        GLOW (t) & 550.8 & 558.8 & 1.5 \\
        PEF (pts) & 3979.6 & 3925.4 & $-$1.4 \\
        \midrule
        \multicolumn{4}{l}{\textit{Design variables:}} \\
        \midrule
        $O/F_{1}$ & 2.9 & 2.9 & 0.6 \\
        $O/F_{2}$ & 3.0 & 3.0 & 0.7 \\
        $m_{\mathrm{prop},1}$ (t) & 397.8 & 425.3 & 6.9 \\
        $m_{\mathrm{prop},2}$ (t) & 89.3 & 72.4 & $-$21.0 \\
        $\dot{m}_{1}$ (kg/s) & 405.5 & 396.7 & $-$2.2 \\
        $\dot{m}_{2}$ (kg/s) & 253.8 & 189.0 & $-$25.5 \\
        $D$ (m) & 4.0 & 4.0 & 0.0 \\
        \midrule
        \multicolumn{4}{l}{\textit{Coupling variables:}} \\
        \midrule
        $I_{sp,vac,1}$ (s) & 339.5 & 339.5 & 0.0 \\
        $I_{sp,vac,2}$ (s) & 362.8 & 362.8 & 0.0 \\
        $A_{e,1}$ (m$^{2}$) & 1.3 & 1.2 & $-$2.2 \\
        $A_{e,2}$ (m$^{2}$) & 4.8 & 3.5 & $-$25.5 \\
        $m_{\mathrm{dry},1}$ (kg) & 33942.1 & 34427.5 & 1.4 \\
        $m_{\mathrm{dry},2}$ (kg) & 6736.8 & 5563.2 & $-$17.4 \\
        1st stage length (m) & 42.0 & 44.7 & 6.4 \\
        2nd stage length (m) & 11.0 & 9.1 & $-$17.1 \\
        GWP (tCO$_2$-eq) & 276.0 & 275.9 & $-$0.1 \\
        Water use (m$^3$world-eq) & 982232.0 & 989204.1 & 0.7 \\
        ADP metals (kgSb-eq) & 16.7 & 16.3 & -2.6 \\
        \bottomrule
    \end{tabular}
\end{table}

Interestingly, most of the PEF values obtained here were not reached in the previous analyses showing that minimizing only the GWP or the water use indicator was not enough for reducing the single score in its current definition, due to their relatively low contribution compared to other indicators such that the ADP for metal and minerals (recall Tab.~\ref{tab:baseline-GWP}). 
In fact, when looking at Tab.~\ref{tab:glow-ss-comparison}, it can be observed that the main driver relies on the reduction of the second stage propellant mass, mass flow rate, and therefore its structural dimensions and its dry mass. 
The climate change and water use impacts remains largely unaffected. These results indicate that, under the current PEF formulation, minimizing the aggregated score is effectively equivalent to minimizing the ADP for metals and minerals, primarily through a reduction of the second-stage mass relative to the first, as explained in the previous subsection. 

\begin{figure*}[!b]
    \centering
\begin{adjustbox}{minipage=\linewidth,scale=1,center}
    \includegraphics[height=7.8cm]{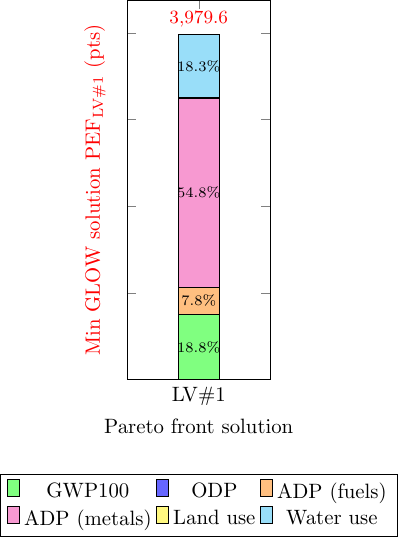}%
    \hfill
    \includegraphics[height=8cm]{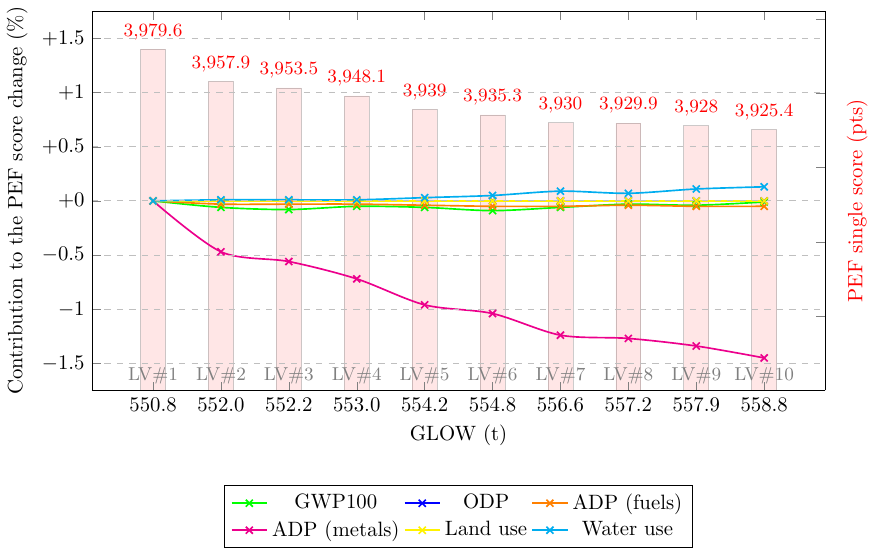}%
\end{adjustbox}
\caption{Left: Contribution of each impact category to the PEF single score for the min GLOW solution ($\mathrm{PEF}_{\mathrm{LV}\#1}$). Right: Evolution of the PEF single score within the multi-objective optimization of GLOW and PEF (background) and contribution of each impact category to the PEF score change with respect to the min GLOW solution (foreground).}
    \label{fig:SS-PEF}
\end{figure*}

This conclusion is observed in Fig.~\ref{fig:SS-PEF} displaying the contributions of each impact category to the PEF score change with respect to the minimum GLOW solution: the contribution of the other indicators to the PEF score change are marginal. 
This behavior can be understood examining Tab.~\ref{tab:PEF_details}, which presents the distribution of PEF and GWP scores by life-cycle stage among the extreme solutions and the median. 
\begin{table}[!h]
    \caption{\\
    Comparison of the PEF and GWP scores distribution for the minimum GLOW, the median, and the minimum PEF solutions. Relative difference $\Delta$ (\%) corresponds to the minimum PEF solution relative to the minimum GLOW solution.
    }
    \label{tab:PEF_details}
    \centering
    \small
    \begin{tabular}{lcccc}
        \toprule
          & Min GLOW & Median & Min PEF & \multirow{2}{*}{$\Delta$ (\%)} \\
          & LV\#1 & LV\#5 & LV\#10 & \\
        \midrule
        PEF (pts) & 3979.6 & 3939.0 & 3925.4 & \\
        \midrule
        Emissions & 243.5 & 242.2 & 243.4 & 0.0 \\
        Propellants  & 973.7 & 981.3 & 991.2 & 1.8 \\
        Components  & 2761.4 & 2714.5 & 2689.7 & $-$2.6 \\
        Transport & 1.1 & 1.1 & 1.1 & 0.0 \\
        \midrule
         \% of total PEF & & & & \\
        \midrule
        Emissions & 6.1 & 6.2 & 6.2 & -- \\
        Propellants & 24.5 & 24.9 & 25.2 & -- \\
        Components & 69.3 & 68.9 & 68.5 & -- \\
        Transport & 0.0 & 0.0 & 0.0 & -- \\
        \midrule
        GWP (tCO\textsubscript{2}-eq) & 276.0 & 275.2 & 275.9 & \\
        \midrule
        Emissions & 88.9 & 88.3 & 88.8 & $-$0.1 \\
        Propellants  & 91.4 & 92.2 & 92.9 & 1.7 \\
        Components  & 95.5 & 94.3 & 93.9 & $-$1.7 \\
        Transport & 0.2 & 0.2 & 0.2 & 0.0 \\
        \bottomrule
    \end{tabular}
\end{table}
In contrast with Tab.~\ref{tab:GWP-details}, where the GWP was minimized alongside the GLOW, launch emissions are no longer the most effective driver for achieving a trade-off between GLOW and PEF score under the current modeling assumptions; instead, the production of launch vehicle components becomes the dominant contributor. This result was expected, since the ADP of metals and minerals contributes the most to the PEF score. Consequently, the optimizer primarily focuses on reducing impacts associated with component production. Interestingly, the contribution of launch emissions remains essentially constant throughout the optimization, while the relative impact of propellant production relatively increases. However, this increase is outweighed by the reduction in impacts from component production, which drives the overall PEF improvement due to its higher contribution. In contrast, the GWP scores remain nearly constant along the optimization, as the gains achieved from lowering component-related impacts are largely offset by the higher impacts associated with propellant production.

Overall, these results show that, with its current definition, minimizing the PEF score alongside the GLOW leads to different design choices compared with the GWP and water-use optimizations studied in the previous sections. While the latter primarily act on propellant use, the PEF optimization places greater emphasis on component production, illustrating how different environmental indicators can drive conflicting design choices in launch vehicle optimization.

\subsection{Synthesis and perspectives of the methodology}

In the previous subsections, the methodology described in Sec.~\ref{sec:3-2} to develop an LCA disciplinary model has been applied to an illustrative MDAO framework for early-stage launch vehicle design, demonstrating the capabilities of integrating environmental considerations. 

For an illustrative test case, multi-objective optimizations were carried out to minimize the launch vehicle GLOW and an environmental impact indicator selected among the GWP, the water use related to the local water scarcity, and the aggregated PEF single score. Overall, the results highlight the antagonistic behaviors between these impact categories. Therefore, it emphasizes the critical importance of carefully selecting the environmental objectives when conducting eco-design studies for launch vehicles, as the choice of indicators can significantly influence the resulting design priorities and trade-offs between different environmental impacts. 
Indeed, while minimizing GWP together with GLOW mainly drives reductions in launch emissions, primarily through adjustments to the oxidizer-to-fuel ratio, minimizing water use leads to the opposite effect. Meanwhile, minimizing the PEF single score, in its current definition, tends to favor configurations with the lowest possible dry masses, at the cost of a higher propellant mass. This analysis reveals a limitation of the single-score approach: because it is strongly dominated by ADP for metals and minerals, it may obscure smaller yet potentially significant variations in other environmental indicators, thereby masking trade-offs at the subsystem level. 
More fundamentally, the analyses show that minimizing a given environmental criterion necessarily requires finding trade-offs between different disciplines, for example between mass-related design choices and trajectory-driven emission estimates. These couplings cannot be consistently captured through a discipline-by-discipline approach and can only be properly addressed within a global MDAO framework integrating an LCA discipline such as the one proposed in this study.
In this context, the latter enables performing multi-objective optimizations across multiple (more than two) environmental indicators (instead of a single score), allowing a more comprehensive analysis of trade-offs among candidate solutions. Furthermore, this enhanced framework also supports other types of early-stage studies, such as sensitivity and parametric analyses \cite{DeOliveira2025}, as well as uncertainty propagation \cite{Brevault2025}.

However, the trends observed in this case study are inherently dependent on its definition. Other formulations could integrate a different launcher architecture (e.g., with 2-3 stages or with boosters, with another propellant mix). In fact, the MDAO architecture could also embed these architectural and technological choices as discrete design variables \cite{Baraton2025, Pelamatti2020}. 
The relatively narrow Pareto fronts obtained in the present study also illustrate the strong constraints imposed by the coupling between vehicle structural mass, propulsion, and trajectory requirements. Although relatively broad ranges were considered for the continuous design variables, the design freedom remains limited when the vehicle architecture and propellant combination are fixed. In particular, allowing the propellant combination itself to be treated as a design variable could introduce additional trade-offs and potentially lead to broader Pareto fronts. More substantial reductions in environmental impacts may therefore require changes beyond the vehicle-level design variables considered here. The MDAO framework integrating an LCA discipline can thus help identify where environmental improvements can be achieved through vehicle design and where changes at a higher architectural or life-cycle level may be required. 

Furthermore, this underlines the dependence of the results observed here on the physical modeling assumptions of the disciplinary models. 
Additionally, the antagonistic behaviors of the environmental indicators highlight the critical importance of carefully defining parametric inventories within the LCA disciplinary model, as subsystem-level modeling assumptions can strongly influence environmental trade-offs. This includes not only the distribution of the components within the launch vehicle architecture and the mass models associated, but also the scope and level of detail of the life-cycle processes that are included (e.g., manufacturing process, transport, type of propellant production). 

As highlighted in Sec.~\ref{sec:2-lca-summary}, significant knowledge gaps remain in accurately assessing the environmental impacts of launch vehicles and the associated uncertainties.  These uncertainties arise not only from the limited availability of space-specific data, but also from several additional sources: the level of fidelity of the LCA discipline, uncertainties in the background LCA databases, uncertainties in the estimation of impact indicators (e.g., GWP definition), and the modeling of launch emissions at high altitude with their corresponding effects, among others.
Among these uncertainties, the environmental effects of high-altitude launch emissions deserve particular attention. Their characterization remains subject to significant knowledge gaps, and recent studies have suggested that their environmental impacts could potentially be substantially larger than those estimated using conventional approaches \cite{DominguezCalabuig2024}. Consequently, the environmental hotspots identified in this study, as well as those obtained from other LCA-based assessments using similar modeling approaches, could shift significantly as the understanding and modeling of high-altitude emissions improve. In particular, the relative contribution of the launch event could become dominant, although future research may also show that its contribution remains limited. This uncertainty should therefore be considered when interpreting the environmental trade-offs and hotspots identified in the present study and highlights the need for improved models and space-specific data on the atmospheric effects of launch emissions.

In addition, the multi-objective optimizations performed in this study illustrate one possible formulation of the MDAO problem under a possible MDAO architecture; identifying the most appropriate strategy requires deeper, case-dependent analyses. 
Nevertheless, the LCA methodology presented here is generic and therefore adaptable to various MDAO problem formulations. Its modular structure also allows the discipline to be progressively improved as more accurate data, emission models, and characterization factors become available. Finally, the parametric life-cycle inventories are adaptable to different launch vehicle architectures and could, in the future, support the integration of additional life-cycle phases, particularly for reusable launch vehicles.

\section{Conclusions}
\label{sec:5}
In this paper, a methodology to integrate a LCA discipline within a MDAO framework dedicated to launch vehicle design has been proposed to enable informed trade-offs between the vehicle performance and its environmental impact. 
It relies on the definition of parametric life-cycle inventories, which depend on the design and coupling variables of the MDAO framework and enable fast computation of the environmental impacts. This approach allows the LCA to be continuously updated along the optimization process. It includes the production of the stages components and propellants, as well as the associated transport to the launch site. Furthermore, launch emissions are evaluated from the optimized trajectory profiles using existing models in the literature \cite{FischerLEATegusphere-2026-3050}, and characterized in terms of climate change impact via GWP, based on available metrics \cite{IPCC2014, Lee2021, Sherwood2018, Hauglustaine2022, Lammel1995}.
The resulting LCA discipline is suitable for MDAO integration and enables environmental considerations to be included at the early design stages of a launch vehicle, alongside other criteria such as overall cost and vehicle performance. The proposed approach has been applied to a representative medium-lift expendable TSTO LOx/LCH\textsubscript{4} launch vehicle design test case illustrating the capabilities of such MDAO framework enhanced with predictive LCA. 

More specifically, multi-objective optimizations were carried out to minimize the launch vehicle GLOW and an environmental impact indicator selected among the GWP, the water use related to the local water scarcity, and the aggregated PEF single score. 
The results reveal the presence of competing effects among impact categories, indicating that environmental optimization cannot be addressed through a single-dimensional perspective. This underscores the need to explicitly define environmental priorities in eco-design studies for launch vehicles, since the selected indicators directly drives design decisions and the resulting trade-offs among environmental impacts.
They also demonstrate the inherent dependency of the analysis on the selected launch vehicle architecture, as well as on the underlying physical modeling assumptions of the disciplines and the level of detail in the life-cycle inventory. 
Nevertheless, the methodology presented in this paper provides the capability to perform multi-criteria optimizations, as well as to incorporate additional design variables (e.g., propellant mix, material distribution), and to progressively improve the disciplinary models as more accurate data become available.
Overall, this paper lays the first foundations for integrating LCA into launch vehicle early-stage design processes, thereby enabling a comprehensive analysis of trade-offs between performance and environmental considerations. 

\section*{Acknowledgements} 
The authors would like to particularly acknowledge Marie Jacquesson, David Miot, Laurence Rozenberg, and Anne-Laure Capomaccio for the valuable insights on life cycle assessment for launch vehicles.

\section*{Declaration of competing interest}
The authors declare that they have no known competing financial interests or personal relationships that could have appeared to influence the work reported in this paper.

\section*{Funding}
The authors acknowledge financial support from the Centre national d’études spatiales (CNES), France (ROR: \url{https://ror.org/04h1h0y33}) through a postdoctoral research grant.

\appendix
\section{Detailed parametric inventory of the representative TSTO LOx/LCH\textsubscript{4} launch vehicle}
\setcounter{table}{0}
\label{app:A}

This section describes the parametric inventories developed to reconstruct the expendable medium-lift TSTO LOx/CH\textsubscript{4} launch vehicle used as the representative test case, following the methodology described in Sec.~\ref{sec:3-2-inventory}. 

The mass of each component is computed in the \textit{Structure} discipline except for the propellant mass which is a design variable. Note that in this study, the fairing mass, $m_\mathrm{fairing}$, is fixed to 2000 kg. For the representative test case, two stages are considered.
The inventory is constructed hierarchically, starting from the individual component inventories and progressively aggregating them into stage-level and launch-vehicle-level inventories.
The first stage accounts for the engine assembly mass (engine, TVC, pressurization system), the thrust frame, the fuel and oxidizer tanks (as well as the intertank if necessary), the TPS associated, and the interstage. Then, the second stage considers the avionics and the EPS modules, the fairing and payload adapter, and all the other components mentioned for the first stage, except the interstage. 

At the component level, each inventory is defined as a parametric foreground inventory combining the relevant background processes for material production and transport. The component inventories are then scaled according to the corresponding component masses and, where applicable, the number of components or user-defined material distribution parameters. At the stage level, these component inventories are aggregated to form the inventory of each stage. The complete launch vehicle inventory is subsequently obtained by aggregating the inventories of the two stages. Propellant inventories are then added separately based on the propellant masses required by each stage.

All the components masses are estimated from the models available in Castellini \cite{castellini2013multidisciplinary}. Note that the thrust frame, interstage, and payload adapter mass models have been modified to account for structures with different material types such as aluminium and composite materials. 
It is also possible to account for a percentage of unused propellant, which is considered as part of the dry mass for the payload performance calculation within the trajectory optimization. In the LCA discipline, however, this propellant is treated separately from the dry mass, and its associated environmental impacts are accounted for as part of the propellant production life-cycle stage.

Using the capabilities of Brightway2 \cite{mutel2017brightway} and lca\_algebraic \cite{jolivet2021lca}, the life-cycle inventory of each launch vehicle component is composed and distributed between the categories \textit{Components} and \textit{Transport}. Tab.~\ref{tab:inventories} summarizes the parametric inventories developed for the representative test case studied in this paper. Note that the material distribution of the engines was established based on a Vulcain-type engine. For the other components, generic inventories were developed based on information available in the literature and are not intended to provide a fully representative LCA of the launch vehicle considered in this study. Further refinement of these inventories through a more detailed LCA of the launch vehicle is left for future work.
\begin{table*}[!p]
    \centering
    \caption{Description of the parametric inventories as defined for this study. 
    They are constructed from the \textit{background} LCI databases: ESA LCA v1.2.0 \cite{ESA_LCA} and ecoinvent v3.9.1 (cutoff) \cite{frischknecht2005ecoinvent}
    }
    \label{tab:inventories}
    \small
    \resizebox{\linewidth}{!}{%
    \begin{tabular}{cccc}
        \toprule
        \multirow{2}{*}{Output from MDO disciplines} & \multirow{2}{*}{Life-cycle category} & \multirow{2}{*}{Distribution parameter} & \multirow{2}{*}{Background process}\\
        \\
        \hline
        \rowcolor{gray!20} \multicolumn{4}{l}{\textit{Common to all stages $i$: tanks, TPS, thrust frame, engine, pressurization system, and propellants}} \\
        \hline
         \multirow{1}{*}{Fuel tank $m_{\mathrm{ftank},i}$} 
            & \multirow{2}{*}{Components} & \multicolumn{1}{c}{--} & \multicolumn{1}{c}{AL2195 generic part, modified from ESA} \\
        \multirow{1}{*}{Oxidizer tank $m_{\mathrm{oxtank},i}$} 
            &                             & \multicolumn{1}{c}{--} & \multicolumn{1}{c}{AL2195 generic part, modified from ESA} \\
        \cline{2-4}
        \multirow{3}{*}{Fuel TPS $m_{\mathrm{TPSf},i}$} 
            & \multirow{3}{*}{Components} 
                & $n_{\mathrm{TPSf},i,\mathrm{CMC}}$ & CMC, ecoinvent\\
            &   & $n_{\mathrm{TPSf},i,\mathrm{foam}}$ & Polyurethane foam, ecoinvent\\
            &   & $n_{\mathrm{TPSf},i,\mathrm{inconel}}$ & Inconel alloy, ESA\\
        \cline{2-4}
        \multirow{2}{*}{Oxidizer TPS $m_{\mathrm{TPSox},i}$} 
            & \multirow{2}{*}{Components} 
                & $n_{\mathrm{TPSox},i,\mathrm{CMC}}$ & CMC, ecoinvent\\
            &   & $n_{\mathrm{TPSox},i,\mathrm{inconel}}$ & Inconel alloy, ESA\\
        \cline{2-4}
        \multirow{2}{*}{Intertank $m_{\mathrm{IT},i}$} 
            & \multirow{2}{*}{Components} 
                & $n_{\mathrm{IT},i,\mathrm{Al}}$ & AL2219 generic part, ESA \\
            &   & $n_{\mathrm{IT},i,\mathrm{CFPR}}$ & CFPR honeycomb 5052, ESA \\
        \hline
        \multirow{2}{*}{$\sum$} 
            & \multirow{2}{*}{Transport} 
                & $d_{\mathrm{tank},i,\mathrm{road}}$ & Lorry 16–32 t EURO4, ecoinvent\\
            & 
                & $d_{\mathrm{tank},i,\mathrm{sea}}$ & Sea container ship, ecoinvent\\
        \hline
         \multirow{4}{*}{Thrust frame $m_{\mathrm{TF},i}$} 
            & \multirow{2}{*}{Components} 
                & $n_{\mathrm{TF},i,\mathrm{Al}}$ & AL2219 generic part, ESA\\
            & 
                & $n_{\mathrm{TF},i,\mathrm{CFPR}}$ & CFPR honeycomb 5052, ESA\\
        \cline{2-4}
            & \multirow{2}{*}{Transport} 
                & $d_{\mathrm{TF},i,\mathrm{road}}$ & Lorry 16–32 t EURO4, ecoinvent \\
            & 
                & $d_{\mathrm{TF},i,\mathrm{sea}}$ & Sea container ship, ecoinvent\\
        \hline
        \multirow{3}{*}{Vulcain-type engines total mass} 
            & \multirow{5}{*}{Components} 
                & $n_{\mathrm{eng},i,\mathrm{Al}}$ & AL2219 generic part, ESA \\
            & 
                & $n_{\mathrm{eng},i,\mathrm{Ni}}$ & Nickel Haynes 25, ESA \\
            & 
                & $n_{\mathrm{eng},i,\mathrm{Ti}}$ & Titanium generic part, ESA \\
        \multirow{2}{*}{$m_{\mathrm{eng},i} \times N_{\mathrm{eng},i}$ } & 
                & $n_{\mathrm{eng},i,\mathrm{FeC}}$ & Stainless steel 304, ESA \\
            & 
                & $n_{\mathrm{eng},i,\mathrm{Cu}}$ & Copper coil, ESA \\
        \cline{2-4}
            & \multirow{2}{*}{Transport} 
                & $d_{\mathrm{eng},i,\mathrm{road}}$ & Lorry 16–32 t EURO4, ecoinvent \\
            & 
                & $d_{\mathrm{eng},i,\mathrm{sea}}$ & Sea container ship, ecoinvent \\
        \hline
        \multirow{1}{*}{Pressurization system mass} 
            & \multirow{3}{*}{Components} & \multicolumn{1}{c}{} \\
        \multirow{1}{*}{Tank $m_{\mathrm{PS,tank},i}$} 
            & & \multicolumn{1}{c}{--} & \multicolumn{1}{c}{AL2195 generic part, modified from ESA} \\
        \multirow{1}{*}{Helium $m_{\mathrm{PS,He},i}$} 
            & & \multicolumn{1}{c}{--} & \multicolumn{1}{c}{Market for Helium, GLO, ecoinvent} \\
        \hline
        \multirow{2}{*}{$\sum$ } 
            & \multirow{2}{*}{Transport} 
                & $d_{\mathrm{PS},i,\mathrm{road}}$ & Lorry 16–32 t EURO4, ecoinvent\\
            & 
                & $d_{\mathrm{PS},i,\mathrm{sea}}$ & Sea container ship, ecoinvent\\
        \hline
        \multirow{1}{*}{Propellant mass} 
            & \multirow{3}{*}{Propellants} & \multicolumn{1}{c}{} \\
        \multirow{1}{*}{Fuel mass $m_{\mathrm{CH_4},i}$}  &  & \multicolumn{1}{c}{--} & CH\textsubscript{4} propellant loading, ESA \\
        \multirow{1}{*}{Oxidizer mass $m_{\mathrm{LOx},i}$} &  & \multicolumn{1}{c}{--} & LOx propellant loading, ESA \\
        \hline
        \rowcolor{gray!20} \multicolumn{4}{l}{\textit{1st stage: interstage included}} \\
        \hline
        \multirow{4}{*}{Interstage mass $m_{\mathrm{IS}}$}
            & \multirow{2}{*}{Components} 
                & $n_{\mathrm{IS},\mathrm{Al}}$ & AL2219 generic part, ESA\\
            & 
                & $n_{\mathrm{IS,CFPR}}$ & CFPR honeycomb 5052, ESA\\
        \cline{2-4}
            & \multirow{2}{*}{Transport} 
                & $d_{\mathrm{IS,road}}$ &Lorry 16–32 t EURO4, ecoinvent\\
            & 
                & $d_\mathrm{{IS,sea}}$ &Sea container ship, ecoinvent\\
        \hline
        \rowcolor{gray!20} \multicolumn{4}{l}{\textit{2nd stage: fairing, payload adapter, avionics and EPS included}} \\
        \hline
            & 
                & $n_{\mathrm{fairing,Al}}$ & AL2219 generic part, ESA \\
            \multirow{-2}{*}{Fairing \& payload adapter} & \multirow{-2}{*}{Components} 
                & $n_{\mathrm{fairing,CFPR}}$ & CFPR honeycomb 5052, ESA \\
        \cline{2-4}
        & 
                & $d_{\mathrm{fairing,road}}$ & Lorry 16–32 t EURO4, ecoinvent \\
        \multirow{-2}{*}{$m_\mathrm{fairing} + m_\mathrm{PLA}$} & \multirow{-2}{*}{Transport} 
                & $d_{\mathrm{fairing,sea}}$ & Sea container ship, ecoinvent\\
        \hline
        \multirow{3}{*}{Avionics mass $m_\mathrm{Av}$} 
            & \multirow{3}{*}{Components} 
                & $n_{\mathrm{Av,Unit}}$ & Electronic unit, ESA \\
            & 
                & $n_{\mathrm{Av,W20}}$ & Single wire AWG20, ESA  \\
            & 
                & $n_{\mathrm{Av,W26}}$ & Single wire AWG26, ESA\\
        \cline{2-4}
        \multirow{3}{*}{ EPS mass $m_\mathrm{EPS}$} 
            & \multirow{3}{*}{Components} 
                & $n_{\mathrm{EPS,Unit}}$ & Power supply unit, ESA \\
            & 
                & $n_{\mathrm{EPS,Bat}}$ & Li-ion battery assembly, ESA\\
            & 
                & $n_{\mathrm{EPS,W14}}$ & Single wire AWG14, ESA \\
        \hline
        \multirow{2}{*}{$\sum$} 
            & \multirow{2}{*}{Transport} 
                & $d_{\mathrm{AvEPS,road}}$ & Lorry 16–32 t EURO4, ecoinvent\\
            & 
                & $d_{\mathrm{AvEPS,sea}}$ & Sea container ship, ecoinvent \\
        \bottomrule
    \end{tabular}
    }
\end{table*}

\textit{Components} accounts for the production of the module from the raw materials selected in the background databases as function of its mass computed in the \textit{Structure} discipline. 
These component-level foreground inventories constitute the elementary building blocks of the stage and launch vehicle inventories described above.
In this study, ESA LCA \cite{ESA_LCA} database is used primarily, but if the material is not found, then the ecoinvent \cite{frischknecht2005ecoinvent} database is used instead. In particular, "modified from ESA" means that since the material was not found, a similar one available in the ESA LCA database is selected and its inventory is modified by replacing the relevant processes using data from the ecoinvent database. When several materials are necessary to reconstruct the component, distribution parameters are allocated with a value chosen by the user. In Eq.~\eqref{eq:eng}, the example is given for the engine mass of the first stage:
\begin{equation}
\label{eq:eng}
    \begin{aligned}
        F_{Components,\mathrm{eng},1} = m_{\mathrm{eng},1} \times \bigg (n_{\mathrm{Al}} \times B_{\mathrm{Al}} + n_{\mathrm{Ni}} \times B_{\mathrm{Ni}}\\
        + n_{\mathrm{Ti}} \times B_{\mathrm{Ti}} + n_{\mathrm{FeC}} \times B_{\mathrm{FeC}} + n_{\mathrm{Cu}} \times B_{\mathrm{Cu}} \bigg)
    \end{aligned}
\end{equation}
where $F_{Components,\mathrm{eng},1}$ denotes the foreground database, $n_{\#},\# = \{\mathrm{Al, Ni, Ti, FeC, Cu}\}$ for aluminum, nickel, titanium, steel, and copper are the distribution parameters and $B_{\#}$ the allocated materials taken from the background databases mentioned. Note that these distribution parameters could also be integrated as design variables within the MDAO architecture. Moreover, the Python library lca\_algebraic \cite{jolivet2021lca} enables the definition of switch parameters to select different raw materials, which could also be embedded in the optimization framework. Similarly, the \textit{Transport} category describes the foreground database which describes the transport phases that are necessary to bring the component to the launch site, as function of its mass and the distance allocated (t.km). In Eq.~\eqref{eq:engtrans}, the example is given for the engine mass of the first stage:
\begin{equation}
\label{eq:engtrans}
    \begin{aligned}
        F_{Transport,\mathrm{eng},1} = m_\mathrm{eng,1} \times \big( d_{\mathrm{road}} \times B_{\mathrm{road}} + d_{\mathrm{sea}} \times B_{\mathrm{sea}} \big)
    \end{aligned}
\end{equation}
where $F_{Transport,\mathrm{eng},1}$ denotes the foreground database, $d_j, j = \{\mathrm{sea}, \mathrm{road}\}$ are the distance traveled by the component for each transport mode, and $B_{j}$ the corresponding background processes. Finally, \textit{Propellants} describes the production of the propellant necessary for the space mission in each stage: it only depends on the design variables related to the propellant mass $m_\mathrm{prop,i}$ and the oxidizer-to-fuel ratio $O/F_i$.

\bibliographystyle{elsarticle-num} 
\bibliography{bibliography}





\end{document}